\documentclass[twoside]{article}
\usepackage{amssymb}
\usepackage{amsmath}

\usepackage[dvips]{epsfig}

\numberwithin{equation}{section}

\begin{document}

\title{{\Large \textbf{\boldmath Waves, damped wave and observation\thanks{%
This work is supported by the NSF of China under grants 10525105,
10771149 and 60974035. Part of this talk was done when the author
visited Fudan University with a financial support from the
"French-Chinese Summer Institute on Applied Mathematics" (September
1-21, 2008).}}}}
\author{{\large Kim Dang PHUNG} \\
%EndAName
{\normalsize \emph{Yangtze Center of Mathematics, Sichuan University,} }\\
{\normalsize \emph{Chengdu 610064, China.} }\\
{\normalsize \emph{E-mail: kim\_dang\_phung@yahoo.fr} }}
\date{\vspace{-12mm}}
\maketitle

\begin{abstract}
This talk describes some applications of two kinds of observation estimate
for the wave equation and for the damped wave equation in a bounded domain
where the geometric control condition of C. Bardos, G. Lebeau and J. Rauch
may failed.
\end{abstract}

%\thispagestyle{first}

% ----------------------------------------------------------------

% ----------------------------------------------------------------
\bigskip

\section{The wave equation and observation}

\bigskip

We consider the wave equation in the solution $u=u(x,t)$%
\begin{equation}
\left\{
\begin{array}{rl}
\partial _{t}^{2}u-\Delta u=0 & \quad \text{in~}\Omega \times \mathbb{R}\
\text{,} \\
u=0 & \quad \text{on~}\partial \Omega \times \mathbb{R}\text{ ,} \\
\left( u,\partial _{t}u\right) \left( \cdot ,0\right) =\left(
u_{0},u_{1}\right) & \text{ ,}%
\end{array}%
\right.  \tag{1.1}  \label{1.1}
\end{equation}%
living in a bounded open set $\Omega $ in $\mathbb{R}^{n}$, $n\geq 1$,
either convex or $C^{2}$\ and connected, with boundary $\partial \Omega $.
It is well-known that for any initial data $\left( u_{0},u_{1}\right) \in
H^{2}(\Omega )\cap H_{0}^{1}\left( \Omega \right) \times H_{0}^{1}\left(
\Omega \right) $, the above problem is well-posed and have a unique strong
solution.

\bigskip

Linked to exact controllability and strong stabilization for the
wave equation (see \cite{Li}), it appears the observability problem
which consists in proving the following estimate%
\[
\left\Vert \left( u_{0},u_{1}\right) \right\Vert _{H_{0}^{1}\left( \Omega
\right) \times L^{2}\left( \Omega \right) }^{2}\leq
C\int_{0}^{T}\int_{\omega }\left\vert \partial _{t}u\left( x,t\right)
\right\vert ^{2}dxdt
\]%
for some constant $C>0$ independent on the initial data. Here, $T>0$ and $%
\omega $\ is a non-empty open subset in $\Omega $. Due to finite
speed of propagation, the time $T$ have to be chosen large enough.
Dealing with high frequency waves i.e., waves which propagates
according the law of geometrical optics, the choice of $\omega $ can
not be arbitrary. In other words, the existence of trapped rays
(e.g, constructed with gaussian beams (see \cite{Ra})) implies the
requirement of some kind of geometric condition on $\left( \omega
,T\right) $ (see \cite{BLR}) in order that the above observability
estimate may hold.

\bigskip

Now, we can ask what kind of estimate we may hope in a geometry with trapped
rays. Let us introduce the quantity
\[
\Lambda =\frac{\left\Vert \left( u_{0},u_{1}\right) \right\Vert _{H^{2}\cap
H_{0}^{1}\left( \Omega \right) \times H_{0}^{1}\left( \Omega \right) }}{%
\left\Vert \left( u_{0},u_{1}\right) \right\Vert _{H_{0}^{1}\left( \Omega
\right) \times L^{2}\left( \Omega \right) }}\text{ ,}
\]%
which can be seen as a measure of the frequency of the wave. In this paper,
we present the two following inequalities%
\begin{equation}
\left\Vert \left( u_{0},u_{1}\right) \right\Vert _{H_{0}^{1}\left( \Omega
\right) \times L^{2}\left( \Omega \right) }^{2}\leq e^{C\Lambda ^{1/\beta
}}\int_{0}^{T}\int_{\omega }\left\vert \partial _{t}u\left( x,t\right)
\right\vert ^{2}dxdt  \tag{1.2}  \label{1.2}
\end{equation}%
and%
\begin{equation}
\left\Vert \left( u_{0},u_{1}\right) \right\Vert _{H_{0}^{1}\left( \Omega
\right) \times L^{2}\left( \Omega \right) }^{2}\leq C\int_{0}^{C\Lambda
^{1/\gamma }}\int_{\omega }\left\vert \partial _{t}u\left( x,t\right)
\right\vert ^{2}dxdt  \tag{1.3}  \label{1.3}
\end{equation}%
where $\beta \in \left( 0,1\right) $, $\gamma >0$. We will also give theirs
applications to control theory.

\bigskip

The strategy to get estimate (\ref{1.2}) is now well-known (see \cite{Ro2},%
\cite{LR}) and a sketch of the proof will be given in Appendix for
completeness. More precisely, we have the following result.

\bigskip

\textbf{Theorem 1.1.-}\qquad \textit{For any }$\omega $\textit{\ non-empty
open subset in }$\Omega $\textit{, for any }$\beta \in \left( 0,1\right) $%
\textit{, there exist }$C>0$\textit{\ and }$T>0$\textit{\ such that for any
solution }$u$\textit{\ of (\ref{1.1}) with non-identically zero initial data
}$\left( u_{0},u_{1}\right) \in H^{2}(\Omega )\cap H_{0}^{1}\left( \Omega
\right) \times H_{0}^{1}\left( \Omega \right) $\textit{, the inequality (\ref%
{1.2}) holds.}

\bigskip

Now, we can ask whether is it possible to get another weight function of $%
\Lambda $ than the exponential one, and in particular a polynomial weight
function with a geometry $\left( \Omega ,\omega \right) $ with trapped rays.
Here we present the following result.

\bigskip

\textbf{Theorem 1.2.-}\qquad \textit{There exists a geometry }$\left( \Omega
,\omega \right) $\textit{\ with trapped rays such that for any solution }$u$%
\textit{\ of (\ref{1.1}) with non-identically zero initial data }$\left(
u_{0},u_{1}\right) \in H^{2}(\Omega )\cap H_{0}^{1}\left( \Omega \right)
\times H_{0}^{1}\left( \Omega \right) $\textit{, the inequality (\ref{1.3})
holds for some }$C>0$\textit{\ and }$\gamma >0$\textit{.}

\bigskip

The proof of Theorem 1.2 is given in \cite{Ph1}. With the help of Theorem
2.1 below, it can also be deduced from \cite{LiR}, \cite{BuH}.

\bigskip

\section{The damped wave equation and our motivation}

\bigskip

We consider the following damped wave equation in the solution $w=w(x,t)$%
\begin{equation}
\left\{
\begin{array}{rl}
\partial _{t}^{2}w-\Delta w+1_{\omega }\partial _{t}w=0 & \quad \text{in~}%
\Omega \times \left( 0,+\infty \right) \ \text{,} \\
w=0 & \quad \text{on~}\partial \Omega \times \left( 0,+\infty \right) \text{
,}%
\end{array}%
\right.  \tag{2.1}  \label{2.1}
\end{equation}%
living in a bounded open set $\Omega $ in $\mathbb{R}^{n}$, $n\geq 1$,
either convex or $C^{2}$\ and connected, with boundary $\partial \Omega $.
Here $\omega $\ is a non-empty open subset in $\Omega $ with trapped rays
and $1_{\omega }$ denotes the characteristic function on $\omega $. Further,
for any $\left( w,\partial _{t}w\right) \left( \cdot ,0\right) \in
H^{2}(\Omega )\cap H_{0}^{1}\left( \Omega \right) \times H_{0}^{1}\left(
\Omega \right) $, the above problem is well-posed for any $t\geq 0$ and have
a unique strong solution.

\bigskip

Denote for any $g\in C\left( \left[ 0,+\infty \right) ;H_{0}^{1}\left(
\Omega \right) \right) \cap C^{1}\left( \left[ 0,+\infty \right)
;L^{2}\left( \Omega \right) \right) $,
\[
E\left( g,t\right) =\frac{1}{2}\int_{\Omega }\left( \left\vert \nabla
g\left( x,t\right) \right\vert ^{2}+\left\vert \partial _{t}g\left(
x,t\right) \right\vert ^{2}\right) dx\text{ .}
\]%
Then for any $0\leq t_{0}<t_{1}$, the strong solution $w$\ satisfies the
following formula%
\begin{equation}
E\left( w,t_{1}\right) -E\left( w,t_{0}\right)
+\int_{t_{0}}^{t_{1}}\int_{\omega }\left\vert \partial _{t}w\left(
x,t\right) \right\vert ^{2}dxdt=0\text{ .}  \tag{2.2}  \label{2.2}
\end{equation}

\bigskip

\subsection{The polynomial decay rate}

\bigskip

Our motivation for establishing estimate (\ref{1.3}) comes from the
following result.

\bigskip

\textbf{Theorem 2.1 .-}\qquad \textit{The following two assertions are
equivalent. Let }$\delta >0$\textit{.}

\begin{description}
\item[\textit{(i)}] \textit{There exists }$C>0$\textit{\ such that for any
solution\ }$w$\textit{\ of (\ref{2.1})\ with the non-null initial data }$%
\left( w,\partial _{t}w\right) \left( \cdot ,0\right) =\left(
w_{0},w_{1}\right) \in H^{2}\left( \Omega \right) \cap H_{0}^{1}\left(
\Omega \right) \times H_{0}^{1}\left( \Omega \right) $\textit{, we have}%
\[
\left\Vert \left( w_{0},w_{1}\right) \right\Vert _{H_{0}^{1}\left( \Omega
\right) \times L^{2}\left( \Omega \right) }^{2}\leq C\int_{0}^{C\left( \frac{%
E\left( \partial _{t}w,0\right) }{E\left( w,0\right) }\right) ^{1/\delta
}}\int_{\omega }\left\vert \partial _{t}w\left( x,t\right) \right\vert
^{2}dxdt\text{ .}
\]

\item[\textit{(ii)}] \textit{There exists }$C>0$\textit{\ such that the
solution }$w$\textit{\ of (\ref{2.1})\ with the initial data }$\left(
w,\partial _{t}w\right) \left( \cdot ,0\right) =\left( w_{0},w_{1}\right)
\in H^{2}\left( \Omega \right) \cap H_{0}^{1}\left( \Omega \right) \times
H_{0}^{1}\left( \Omega \right) $\textit{\ satisfies}%
\[
E\left( w,t\right) \leq \frac{C}{t^{\delta }}~\left\Vert \left(
w_{0},w_{1}\right) \right\Vert _{H^{2}\cap H_{0}^{1}\left( \Omega \right)
\times H_{0}^{1}\left( \Omega \right) }^{2}\quad \forall t>0\text{ .}
\]
\end{description}

\bigskip

\textbf{Remark .-}\qquad It is not difficult to see (e.g., \cite{Ph2}) by a
classical decomposition method, a translation in time and (\ref{2.2}), that
the inequality (\ref{1.3}) with the exponent $\gamma $ for the wave equation
implies the inequality of $(i)$ in Theorem 2.1 with the exponent $\delta
=2\gamma /3$ for the damped wave equation. And conversely, the inequality of
$(i)$ in Theorem 2.1 with the exponent $\delta $ for the damped wave
equation implies the inequality (\ref{1.3}) with the exponent $\gamma
=\delta /2$ for the wave equation.

\bigskip

Proof of Theorem 2.1.-

$\left( ii\right) \Rightarrow \left( i\right) $. Suppose that
\[
E\left( w,T\right) \leq \frac{C}{T^{\delta }}\left\Vert \left(
w_{0},w_{1}\right) \right\Vert _{H^{2}\cap H_{0}^{1}\left( \Omega \right)
\times H_{0}^{1}\left( \Omega \right) }^{2}\quad \forall T>0\text{ .}
\]%
Therefore from (\ref{2.2})
\[
E\left( w,0\right) \leq \frac{C}{T^{\delta }}\left\Vert \left(
w_{0},w_{1}\right) \right\Vert _{H^{2}\cap H_{0}^{1}\left( \Omega \right)
\times H_{0}^{1}\left( \Omega \right) }^{2}+\int_{0}^{T}\int_{\omega
}\left\vert \partial _{t}w\left( x,t\right) \right\vert ^{2}dxdt\text{ .}
\]%
By choosing%
\[
T=\left( 2C\frac{\left\Vert \left( w_{0},w_{1}\right) \right\Vert
_{H^{2}\cap H_{0}^{1}\left( \Omega \right) \times H_{0}^{1}\left( \Omega
\right) }^{2}}{E\left( w,0\right) }\right) ^{1/\delta }\text{ ,}
\]%
we get the desired estimate
\[
E\left( w,0\right) \leq 2\int_{0}^{\left[ 2C\frac{\left\Vert \left(
w_{0},w_{1}\right) \right\Vert _{H^{2}\cap H_{0}^{1}\left( \Omega \right)
\times H_{0}^{1}\left( \Omega \right) }^{2}}{E\left( w,0\right) }\right]
^{1/\delta }}\int_{\omega }\left\vert \partial _{t}w\left( x,t\right)
\right\vert ^{2}dxdt\text{ .}
\]

$\left( i\right) \Rightarrow \left( ii\right) $. Conversely, suppose the
existence of a constant $c>1$ such that the solution $w$ of (\ref{2.1}) with
the non-null initial data $\left( w,\partial _{t}w\right) \left( \cdot
,0\right) =\left( w_{0},w_{1}\right) \in H^{2}\left( \Omega \right) \cap
H_{0}^{1}\left( \Omega \right) \times H_{0}^{1}\left( \Omega \right) $\
satisfies
\[
E\left( w,0\right) \leq c\int_{0}^{c\left( \frac{E\left( w,0\right) +E\left(
\partial _{t}w,0\right) }{E\left( w,0\right) }\right) ^{1/\delta
}}\int_{\omega }\left\vert \partial _{t}w\left( x,t\right) \right\vert
^{2}dxdt\text{ .}
\]%
We obtain the following inequalities by a translation on the time variable
and by using (\ref{2.2}). $\forall s\geq 0$
\[
\begin{array}{ll}
\frac{E\left( w,s\right) }{\left( E\left( w,0\right) +E\left( \partial
_{t}w,0\right) \right) } & \leq c\int_{s}^{s+c\left( \frac{\left( E\left(
w,0\right) +E\left( \partial _{t}w,0\right) \right) }{E\left( w,s\right) }%
\right) ^{1/\delta }}\int_{\omega }\frac{\left\vert \partial _{t}w\left(
x,t\right) \right\vert ^{2}}{E\left( w,0\right) +E\left( \partial
_{t}w,0\right) }dxdt \\
& \leq c\left( \frac{E\left( w,s\right) }{E\left( w,0\right) +E\left(
\partial _{t}w,0\right) }-\frac{E\left( w,s+c\left( \frac{E\left( w,0\right)
+E\left( \partial _{t}w,0\right) }{E\left( w,s\right) }\right) ^{1/\delta
}\right) }{E\left( w,0\right) +E\left( \partial _{t}w,0\right) }\right)
\text{ .}%
\end{array}%
\]%
Denoting $G\left( s\right) =\frac{E\left( w,s\right) }{E\left( w,0\right)
+E\left( \partial _{t}w,0\right) }$, we deduce using the decreasing of $G$
that%
\[
G\left( s+c\left( \frac{1}{G\left( s\right) }\right) ^{1/\delta }\right)
\leq G\left( s\right) \leq c\left[ G\left( s\right) -G\left( s+c\left( \frac{%
1}{G\left( s\right) }\right) ^{1/\delta }\right) \right]
\]%
which gives%
\[
G\left( s+c\left( \frac{1}{G\left( s\right) }\right) ^{1/\delta }\right)
\leq \frac{c}{1+c}G\left( s\right) \text{ .}
\]%
Let $c_{1}=\left( \frac{1+c}{c}\right) ^{1/\delta }-1>0$ and denoting $%
d\left( s\right) =\left( \frac{c}{c_{1}}\frac{1}{s}\right) ^{\delta }$. We
distinguish two cases.

\noindent If $c_{1}s\leq c\left( \frac{1}{G\left( s\right) }\right)
^{1/\delta }$, then $G\left( s\right) \leq \left( \frac{c}{c_{1}}\frac{1}{s}%
\right) ^{\delta }$ and%
\[
G\left( \left( 1+c_{1}\right) s\right) \leq d\left( s\right) \text{ .}
\]

\noindent If $c_{1}s>c\left( \frac{1}{G\left( s\right) }\right) ^{1/\delta }$%
, then $s+c\left( \frac{1}{G\left( s\right) }\right) ^{1/\delta }<\left(
1+c_{1}\right) s$ and the decreasing of $G$ gives $G\left( \left(
1+c_{1}\right) s\right) \leq G\left( s+c\left( \frac{1}{G\left( s\right) }%
\right) ^{1/\delta }\right) $ and then%
\[
G\left( \left( 1+c_{1}\right) s\right) \leq \frac{c}{1+c}G\left( s\right)
\text{ .}
\]%
Consequently, we have that $\forall s>0$, $\forall n\in \mathbb{N}$, $n\geq
1 $,%
\[
\begin{array}{ll}
G\left( \left( 1+c_{1}\right) s\right) & \leq \max \left[ d\left( s\right) ,%
\frac{c}{1+c}d\left( \frac{s}{\left( 1+c_{1}\right) }\right) ,\cdot \cdot
\cdot ,\right. \\
& \qquad \qquad \left. ,\left( \frac{c}{1+c}\right) ^{n}d\left( \frac{s}{%
\left( 1+c_{1}\right) ^{n}}\right) ,\left( \frac{c}{1+c}\right)
^{n+1}G\left( \frac{s}{\left( 1+c_{1}\right) ^{n}}\right) \right] \text{ .}%
\end{array}%
\]%
Now, remark that with our choice of $c_{1}$, we get%
\[
\frac{c}{1+c}d\left( \frac{s}{\left( 1+c_{1}\right) }\right) =d\left(
s\right) \quad \forall s>0\text{ .}
\]%
Thus, we deduce that $\forall n\geq 1$
\[
\begin{array}{ll}
G\left( \left( 1+c_{1}\right) s\right) & \leq \max \left( d\left( s\right)
,\left( \frac{c}{1+c}\right) ^{n+1}G\left( \frac{s}{\left( 1+c_{1}\right)
^{n}}\right) \right) \\
& \leq \max \left( d\left( s\right) ,\left( \frac{c}{1+c}\right)
^{n+1}\right) \qquad \text{because }G\leq 1\text{ ,}%
\end{array}%
\]%
and conclude that $\forall s>0$
\[
\frac{E\left( w,s\right) }{E\left( w,0\right) +E\left( \partial
_{t}w,0\right) }=G\left( s\right) \leq d\left( \frac{s}{1+c_{1}}\right)
=\left( \frac{c\left( 1+c_{1}\right) }{c_{1}}\right) ^{\delta }\frac{1}{%
s^{\delta }}\text{ .}
\]%
This completes the proof.

\bigskip

\subsection{The approximate controllability}

\bigskip

The goal of this section consists in giving an application of estimate (\ref%
{1.2}).

\bigskip

\noindent For any $\omega $\ non-empty open subset in $\Omega $, for any $%
\beta \in \left( 0,1\right) $, let $T>0$ be given in Theorem 1.1.

\bigskip

\noindent Let $(v_{0},v_{1},v_{0d},v_{1d})\in \left( H^{2}(\Omega
)\cap H_{0}^{1}\left( \Omega \right) \times H_{0}^{1}\left( \Omega
\right) \right) ^{2}$ and $u$ be the solution of (\ref{1.1}) with
initial data $\left( u,\partial _{t}u\right) (\cdot
,0)=(v_{0},v_{1})$.

\bigskip

\noindent For any integer $N>0$, let us introduce
\begin{equation}
f_{N}\left( x,t\right) =-1_{\omega }\sum\limits_{\ell =0}^{N}\left[ \partial
_{t}w^{\left( 2\ell +1\right) }\left( x,t\right) +\partial _{t}w^{\left(
2\ell \right) }\left( x,T-t\right) \right] \text{ ,}  \tag{2.3}  \label{2.3}
\end{equation}%
where $w^{(0)}\in C\left( \left[ 0,T\right] ;H^{2}(\Omega )\cap
H_{0}^{1}\left( \Omega \right) \right) $ is the solution of the damped wave
equation (\ref{2.1}) with initial data
\[
\left( w^{\left( 0\right) },\partial _{t}w^{\left( 0\right) }\right) \left(
\cdot ,0\right) =(v_{0d},-v_{1d})-\left( u,-\partial _{t}u\right) (\cdot ,T)%
\text{ in }\Omega \text{ ,}
\]%
and for $j\geq 0$, $w^{(j+1)}\in C\left( \left[ 0,T\right] ;H^{2}(\Omega
)\cap H_{0}^{1}\left( \Omega \right) \right) $ is the solution of the damped
wave equation (\ref{2.1}) with initial data
\[
\left( w^{\left( j+1\right) },\partial _{t}w^{\left( j+1\right) }\right)
\left( \cdot ,0\right) =\left( -w^{\left( j\right) },\partial _{t}w^{\left(
j\right) }\right) \left( \cdot ,T\right) \text{ in }\Omega \text{ .}
\]%
Introduce
\[
M=\underset{j\geq 0}{\sup }\left\Vert w^{\left( j\right) }\left( \cdot
,0\right) ,\partial _{t}w^{\left( j\right) }\left( \cdot ,0\right)
\right\Vert _{H^{2}\left( \Omega \right) \times H_{0}^{1}(\Omega )}^{2}\text{
.}
\]

\bigskip

\noindent Our main result is as follows.

\bigskip

\textbf{Theorem 2.2 .-}\qquad \textit{Suppose that }$M<+\infty $\textit{.
Then there exists }$C>0$\textit{\ such that for all }$N>0$\textit{, the
control function }$f_{N}$\textit{\ given by (\ref{2.3}), drives the system }%
\[
\left\{
\begin{array}{ll}
\partial _{t}^{2}v-\Delta v=1_{\omega \times \left( 0,T\right) }f_{N} &
\quad \text{\textit{in} }\Omega \times \left( 0,T\right) \text{ ,} \\
v=0 & \quad \text{\textit{on} }\partial \Omega \times \left( 0,T\right)
\text{ ,} \\
\left( v,\partial _{t}v\right) (\cdot ,0)=(v_{0},v_{1}) & \quad \text{%
\textit{in} }\Omega \text{ ,}%
\end{array}%
\right.
\]%
\textit{to the desired data }$(v_{0d},v_{1d})$\textit{\ approximately at
time }$T$\textit{\ i.e., }%
\[
\left\Vert v\left( \cdot ,T\right) -v_{0d},\partial _{t}v\left( \cdot
,T\right) -v_{1d}\right\Vert _{H_{0}^{1}\left( \Omega \right) \times
L^{2}(\Omega )}^{2}\leq \frac{C}{\left[ \ln \left( 1+2N\right) \right]
^{2\beta }}M\text{ ,}
\]%
\textit{and satisfies }%
\[
\left\Vert f_{N}\right\Vert _{L^{\infty }\left( 0,T;L^{2}\left( \Omega
\right) \right) }\leq C\left( N+1\right) \left\Vert
(v_{0},v_{1},v_{0d},v_{1d})\right\Vert _{\left( H_{0}^{1}\left( \Omega
\right) \times L^{2}(\Omega )\right) ^{2}}\text{ .}
\]

\bigskip

\textbf{Remark .-}\qquad For any $\varepsilon >0$, we can choose $N$ such
that
\[
\frac{C}{\left[ \ln \left( 1+2N\right) \right] ^{2\beta }}M\simeq
\varepsilon ^{2}\text{ and }\left( 2N+1\right) \simeq e^{\left( \frac{\sqrt{%
CM}}{\varepsilon }\right) ^{1/\beta }}\text{ ,}
\]%
in order that
\[
\left\Vert v\left( \cdot ,T\right) -v_{0d},\partial _{t}v\left( \cdot
,T\right) -v_{1d}\right\Vert _{H_{0}^{1}\left( \Omega \right) \times
L^{2}(\Omega )}\leq \varepsilon \text{ ,}
\]%
and
\[
\left\Vert f\right\Vert _{L^{\infty }\left( 0,T;L^{2}\left( \Omega \right)
\right) }\leq e^{\left[ \left( \frac{C}{\varepsilon }\sqrt{M}\right)
^{1/\beta }\right] }\left\Vert (v_{0},v_{1},v_{0d},v_{1d})\right\Vert
_{\left( H_{0}^{1}\left( \Omega \right) \times L^{2}(\Omega )\right) ^{2}}%
\text{ .}
\]%
In \cite{Zu}, a method was proposed to construct an approximate control. It
consists of minimizing a functional depending on the parameter $\varepsilon $%
. However there, no estimate of the cost is given. On the other
hand, estimate of the form (\ref{1.2}) was originally established by
\cite{Ro2} to give the cost (see also \cite{Le}). Here, we present a
new way to construct an approximate control by superposing different
waves. Given a cost to not overcome, we construct a solution which
will be closed in the above sense to the desired state. It takes
ideas from \cite{Ru} and \cite{BF} like an iterative time reversal
construction (see also \cite{PPV} and \cite{ZL}).\bigskip

\bigskip

\subsubsection{Proof}

\bigskip

Consider the solution
\[
V\left( \cdot ,t\right) =\sum\limits_{\ell =0}^{N}\left[ w^{\left( 2\ell
+1\right) }\left( \cdot ,t\right) +w^{\left( 2\ell \right) }\left( \cdot
,T-t\right) \right] \text{ .}
\]%
We deduce that for $t\in \left( 0,T\right) $%
\[
\left\{
\begin{array}{ll}
\partial _{t}^{2}V\left( \cdot ,t\right) -\Delta V\left( \cdot ,t\right)
=-1_{\omega }\sum\limits_{\ell =0}^{N}\left[ \partial _{t}w^{\left( 2\ell
+1\right) }\left( \cdot ,t\right) +\partial _{t}w^{\left( 2\ell \right)
}\left( \cdot ,T-t\right) \right] \text{ ,} &  \\
V=0\quad \text{on }\partial \Omega \times \left( 0,T\right) \text{ ,} &  \\
\left( V,\partial _{t}V\right) \left( \cdot ,0\right) =0\quad \text{in }%
\Omega \text{\textit{\ .}} &
\end{array}%
\right.
\]%
Now, from the definition of $w^{\left( 0\right) }$, the property of $\left(
w^{\left( j+1\right) },\partial _{t}w^{\left( j+1\right) }\right) \left(
\cdot ,0\right) $ and a change of variable, we obtain that%
\[
\begin{array}{ll}
\left( V,\partial _{t}V\right) \left( \cdot ,T\right) & =\left( w^{\left(
0\right) },-\partial _{t}w^{\left( 0\right) }\right) \left( \cdot ,0\right)
+\left( w^{\left( 2N+1\right) },\partial _{t}w^{\left( 2N+1\right) }\right)
\left( \cdot ,T\right) \\
& =(v_{0d},v_{1d})-\left( u,\partial _{t}u\right) (\cdot ,T)+\left(
w^{\left( 2N+1\right) },\partial _{t}w^{\left( 2N+1\right) }\right) \left(
\cdot ,T\right) \text{ .}%
\end{array}%
\]%
Finally, the solution $v=V+u$ satisfies
\[
\left\{
\begin{array}{ll}
\partial _{t}^{2}v-\Delta v=1_{\omega \times \left( 0,T\right) }f_{N}\quad
\text{in }\Omega \times \left( 0,T\right) \text{ ,} &  \\
v=0\quad \text{on }\partial \Omega \times \left( 0,T\right) \text{ ,} &  \\
\left( v,\partial _{t}v\right) \left( \cdot ,0\right) =(v_{0},v_{1})\quad
\text{in }\Omega \text{ ,} &  \\
\left( v,\partial _{t}v\right) \left( \cdot ,T\right)
=(v_{0d},v_{1d})+\left( w^{\left( 2N+1\right) },\partial _{t}w^{\left(
2N+1\right) }\right) \left( \cdot ,T\right) \quad \text{in }\Omega \text{%
\textit{\ .}} &
\end{array}%
\right.
\]%
Clearly,
\[
\left\Vert v\left( \cdot ,T\right) -v_{0d},\partial _{t}v\left( \cdot
,T\right) -v_{1d}\right\Vert _{H_{0}^{1}\left( \Omega \right) \times
L^{2}(\Omega )}^{2}=2E\left( w^{\left( 2N+1\right) },T\right) \text{ .}
\]%
It remains to estimate $E\left( w^{\left( 2N+1\right) },T\right) $. We claim
that
\[
\exists C>0\qquad \forall N\geq 1\qquad E\left( w^{\left( 2N+1\right)
},T\right) \leq \frac{C}{\left[ \ln \left( 1+2N\right) \right] ^{2\beta }}M%
\text{ .}
\]

\bigskip

Indeed, from Theorem 1.1, we can easily see by a classical decomposition
method that there exist $C>0$ and $T>0$ such that for any $j\geq 0$,
\[
\begin{array}{ll}
& \quad \left\Vert w^{\left( j+1\right) }\left( \cdot ,0\right) ,\partial
_{t}w^{\left( j+1\right) }\left( \cdot ,0\right) \right\Vert
_{H_{0}^{1}\left( \Omega \right) \times L^{2}\left( \Omega \right) }^{2} \\
& \leq C\exp \left( C\frac{\left\Vert w^{\left( j+1\right) }\left( \cdot
,0\right) ,\partial _{t}w^{\left( j+1\right) }\left( \cdot ,0\right)
\right\Vert _{H^{2}\left( \Omega \right) \times H^{1}\left( \Omega \right) }%
}{\left\Vert w^{\left( j+1\right) }\left( \cdot ,0\right) ,\partial
_{t}w^{\left( j+1\right) }\left( \cdot ,0\right) \right\Vert
_{H_{0}^{1}\left( \Omega \right) \times L^{2}\left( \Omega \right) }}\right)
^{1/\beta } \\
& \quad \int_{0}^{T}\int_{\omega }\left\vert \partial _{t}w^{\left(
j+1\right) }\left( x,t\right) \right\vert ^{2}dxdt\text{ .}%
\end{array}%
\]%
Since
\[
E\left( w^{\left( j+1\right) },0\right) =E\left( w^{\left( j\right)
},T\right) \qquad \forall j\geq 0\text{ ,}
\]%
we deduce from (\ref{2.2}) that for any $j\geq 0$%
\[
\begin{array}{ll}
E\left( w^{\left( j+1\right) },0\right) & \leq C\exp \left( C\frac{M}{%
\left\Vert w^{\left( j+1\right) }\left( \cdot ,0\right) ,\partial
_{t}w^{\left( j+1\right) }\left( \cdot ,0\right) \right\Vert
_{H_{0}^{1}\left( \Omega \right) \times L^{2}\left( \Omega \right) }^{2}}%
\right) ^{1/\left( 2\beta \right) } \\
& \quad \left[ E\left( w^{\left( j\right) },T\right) -E\left( w^{\left(
j+1\right) },T\right) \right] \text{ .}%
\end{array}%
\]%
Let
\[
d_{j}=E\left( w^{\left( j+1\right) },T\right) \text{ .}
\]%
By using the decreasing property of the sequence $d_{j}$, that is $d_{j}\leq
d_{j-1}$, we obtain that for any integer $0\leq j\leq 2N$
\[
d_{j}\leq Ce^{\left( C\frac{M}{d_{2N}}\right) ^{1/\left( 2\beta \right) }}%
\left[ d_{j-1}-d_{j}\right] \text{ .}
\]%
By summing over $\left[ 0,2N\right] $, we deduce that
\[
\left( 2N+1\right) d_{2N}\leq Ce^{\left( C\frac{M}{d_{2N}}\right) ^{1/\left(
2\beta \right) }}\left[ d_{-1}-d_{2N}\right] \text{ .}
\]%
Finally, using the fact that $d_{-1}\leq M$, it follows that
\[
d_{2N}\leq \frac{C}{\left[ \ln \left( 1+2N\right) \right] ^{2\beta }}M\text{
.}
\]%
This completes the proof of our claim.

\bigskip

On the other hand, the computation of the bound of $f_{N}$ is immediate.
Therefore, we check that for some $C>0$ and $T>0$,%
\[
\left\Vert f_{N}\right\Vert _{L^{\infty }\left( 0,T;L^{2}\left( \Omega
\right) \right) }\leq C\left( N+1\right) \left\Vert
(v_{0},v_{1},v_{0d},v_{1d})\right\Vert _{\left( H_{0}^{1}\left( \Omega
\right) \times L^{2}(\Omega )\right) ^{2}}\text{ ,}
\]%
\[
\left\Vert v\left( \cdot ,T\right) -v_{0d},\partial _{t}v\left( \cdot
,T\right) -v_{1d}\right\Vert _{H_{0}^{1}\left( \Omega \right) \times
L^{2}(\Omega )}^{2}\leq \frac{C}{\left[ \ln \left( 1+2N\right) \right]
^{2\beta }}M\text{ ,}
\]%
for any $\beta \in \left( 0,1\right) $ and any integer $N>0$. This completes
the proof of our Theorem.

\bigskip

\subsubsection{Numerical experiments}

\bigskip

Here, we perform numerical experiments to investigate the practical
applicability of the approach proposed to construct an approximate control.
For simplicity, we consider a square domain $\Omega =\left( 0,1\right)
\times \left( 0,1\right) $, $\omega =\left( 0,1/5\right) \times \left(
0,1\right) $. The time of controllability is given by $T=4$.

\bigskip

\noindent For convenience we recall some well-known formulas. Denote by $%
\left\{ e_{j}\right\} _{j\geq1}$ the Hilbert basis in $L^{2}\left(
\Omega\right) $ formed by the eigenfunctions of the operator $-\Delta$ with
eigenvalues $\left\{ \lambda_{j}\right\} _{j\geq1}$, such that $\left\|
e_{j}\right\| _{L^{2}\left( \Omega\right) }=1$ and $0<\lambda_{1}<%
\lambda_{2}\leq\lambda_{3}\leq\cdots$, i.e.,
\[
\left\{
\begin{array}{ll}
\lambda_{j}=\pi^{2}\left( k_{j}^{2}+\ell_{j}^{2}\right) \text{, }\qquad
k_{j},\ell_{j}\in\mathbb{N}^{\ast}\text{,} &  \\
e_{j}\left( x_{1},x_{2}\right) =2\sin\left( \pi k_{j}x_{1}\right) \sin\left(
\pi\ell_{j}x_{2}\right) \text{ .} &
\end{array}
\right.
\]

\noindent The solution of
\[
\left\{
\begin{array}{rl}
\partial _{t}^{2}v-\Delta v=f & \quad \text{in }\Omega \times \left(
0,T\right) \text{ ,} \\
v=0 & \quad \text{on }\partial \Omega \times \left( 0,T\right) \text{ ,} \\
\left( v,\partial _{t}v\right) (\cdot ,0)=(v_{0},v_{1}) & \quad \text{in }%
\Omega \text{ ,}%
\end{array}%
\right.
\]

\noindent where $f$ is in the form
\[
f\left( x_{1},x_{2}\right) =-1_{\omega }\sum\limits_{j\geq 1}f_{j}\left(
t\right) e_{j}\left( x_{1},x_{2}\right) \text{ ,}
\]

\noindent is given by the formula%
\[
\begin{array}{cc}
& v\left( x_{1},x_{2},t\right) =\underset{G\rightarrow +\infty }{\lim }%
\sum\limits_{j=1}^{G}\left\{ a_{j}^{0}\cos \left( t\sqrt{\lambda _{j}}%
\right) +a_{j}^{1}\frac{1}{\sqrt{\lambda _{j}}}\sin \left( t\sqrt{\lambda
_{j}}\right) \quad \quad \right. \\
\multicolumn{1}{r}{} & \multicolumn{1}{r}{\left. +\frac{1}{\sqrt{\lambda _{j}%
}}\int_{0}^{t}\sin \left( \left( t-s\right) \sqrt{\lambda _{j}}\right)
R_{j}\left( s\right) ds\right\} e_{j}\left( x_{1},x_{2}\right) \text{ ,}}%
\end{array}%
\]

\noindent where
\[
\left\{
\begin{array}{l}
v_{0}\left( x_{1},x_{2}\right) =\underset{G\rightarrow +\infty }{\lim }%
\sum\limits_{j=1}^{G}a_{j}^{0}e_{j}\left( x_{1},x_{2}\right) \text{,}\quad
\sum\limits_{j\geq 1}\lambda _{j}\left\vert a_{j}^{0}\right\vert
^{2}<+\infty \text{ ,} \\
v_{1}\left( x_{1},x_{2}\right) =\underset{G\rightarrow +\infty }{\lim }%
\sum\limits_{j=1}^{G}a_{j}^{1}e_{j}\left( x_{1},x_{2}\right) \text{,}\quad
\sum\limits_{j\geq 1}\left\vert a_{j}^{1}\right\vert ^{2}<+\infty \text{ ,}
\\
R_{j}\left( t\right) =-\underset{G\rightarrow +\infty }{\lim }%
\sum\limits_{i=1}^{G}\left( \int_{\omega }e_{i}e_{j}dx_{1}dx_{2}\right)
f_{i}\left( t\right) \text{ .}%
\end{array}%
\right.
\]

\noindent Here, $G$ will be the number of Galerkin mode. The numerical
results are shown below. The approximate solution of the damped wave
equation is established via a system of ODE solved by MATLAB.

\bigskip

\paragraph{Example 1 : low frequency.}

\bigskip

The initial condition and desired target are specifically as follows. $%
\left( v_{0},v_{1}\right) =\left( 0,0\right) $ and $\left(
v_{0d},v_{1d}\right) =\left( e_{1}+e_{2},e_{1}\right) $. We take the number
of Galerkin mode $G=100$ and the number of iterations in the time reversal
construction $N=30$.

\bigskip

Below, we plot the graph of the desired initial data $v_{0d}$ and the
controlled solution $v\left( \cdot,t=T=4\right) $.

\[
\begin{array}{cc}
\raisebox{-0cm}{\includegraphics[ height=4.0918cm, width=5.3065cm
]{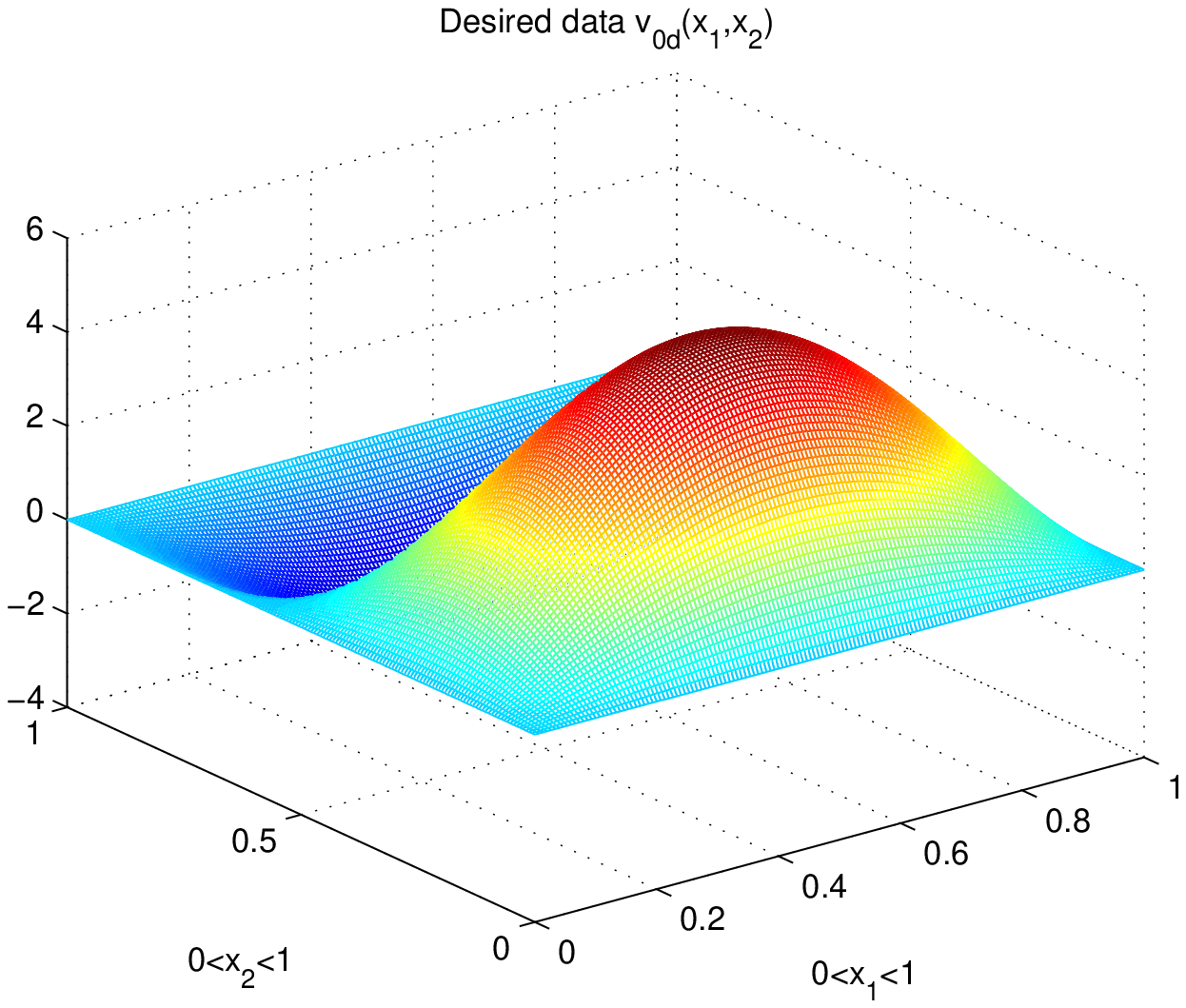}} &
\raisebox{-0cm}{\includegraphics[ height=4.0918cm, width=5.3065cm
]{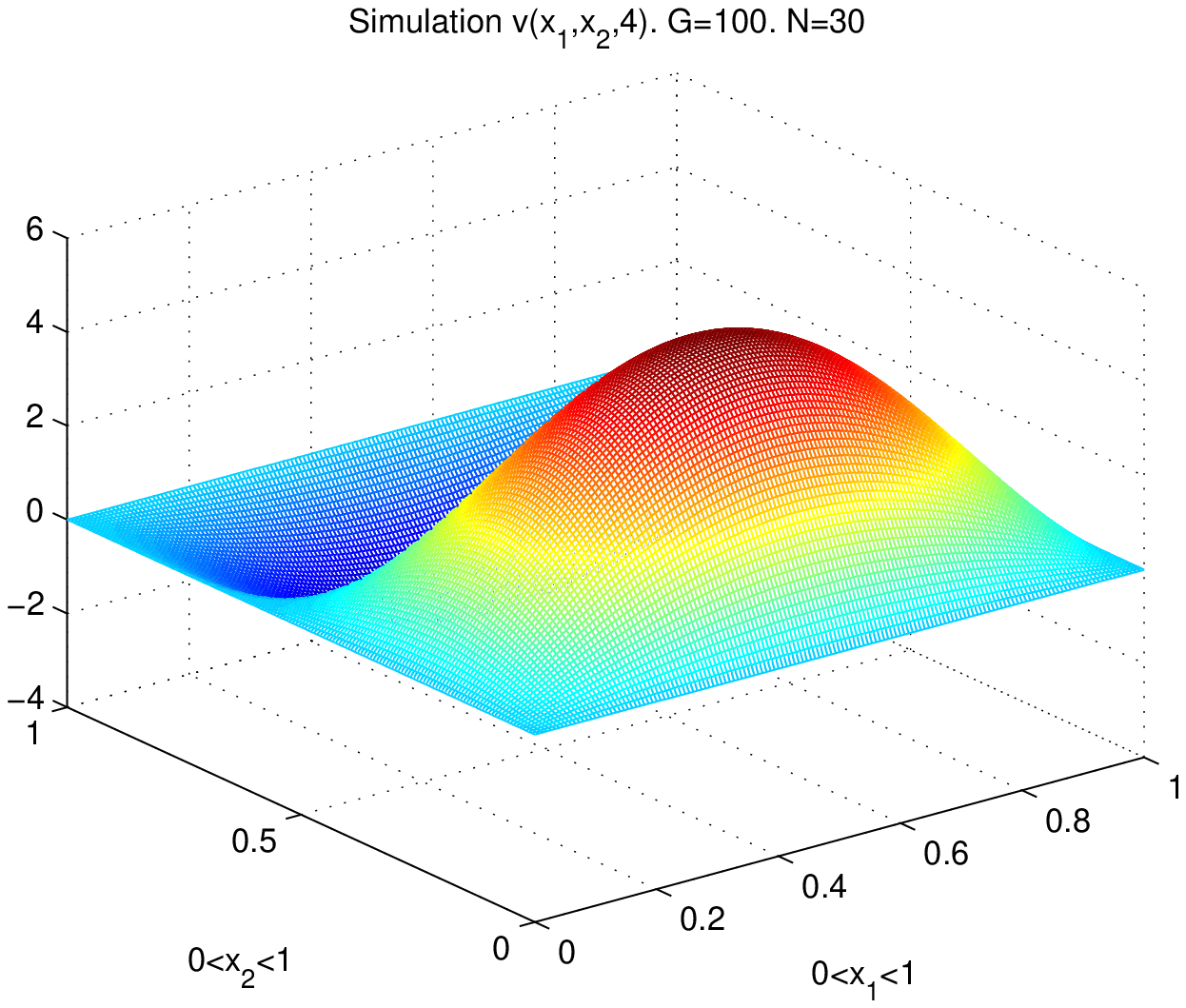}}%
\end{array}
\]

Below, we plot the graph of the energy of the controlled solution and the
cost of the control function.

\[
\begin{array}{cc}
\raisebox{-0cm}{\includegraphics[ height=4.0918cm, width=5.3065cm
]{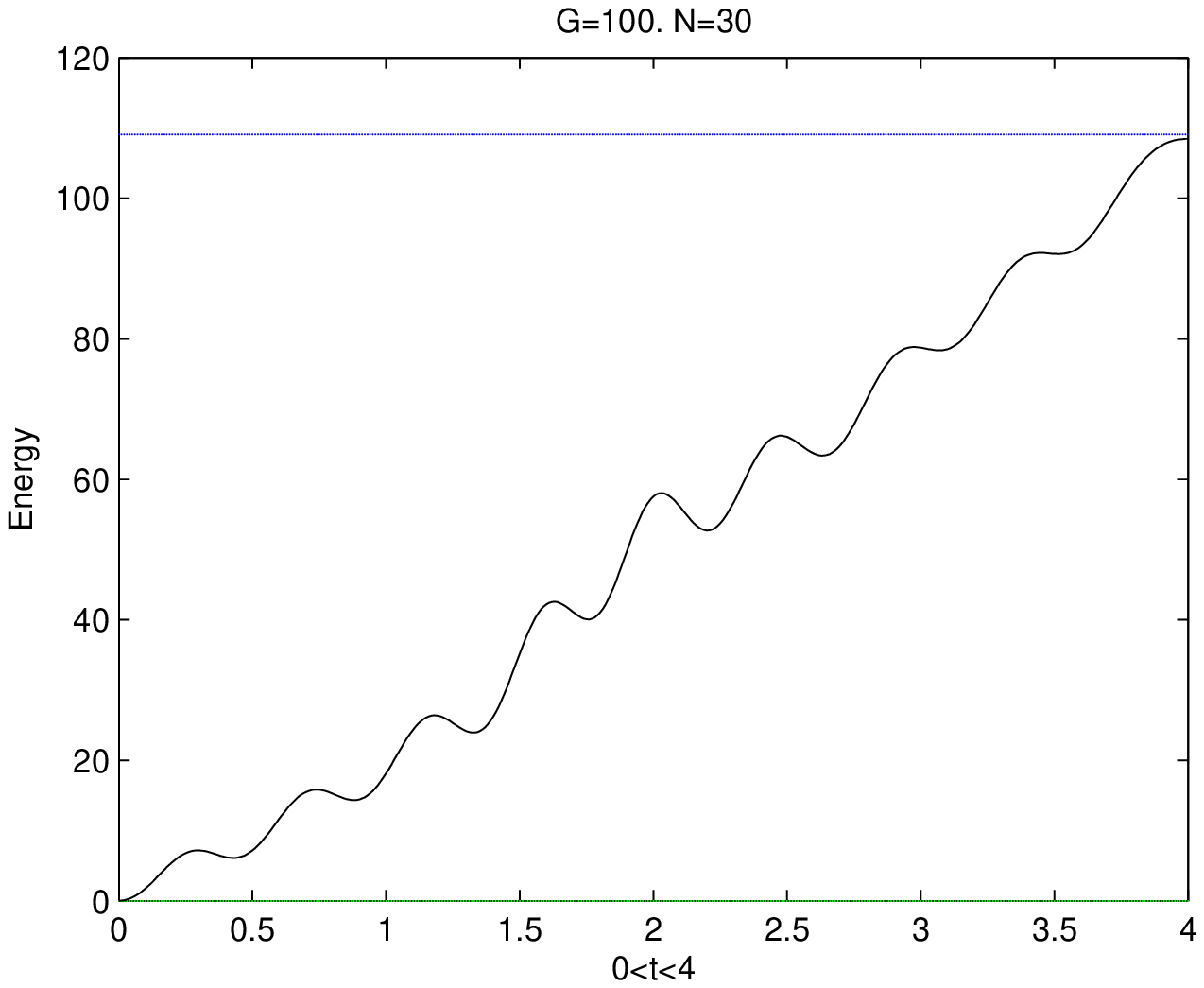}} &
\raisebox{-0cm}{\includegraphics[ height=4.0918cm, width=5.3065cm
]{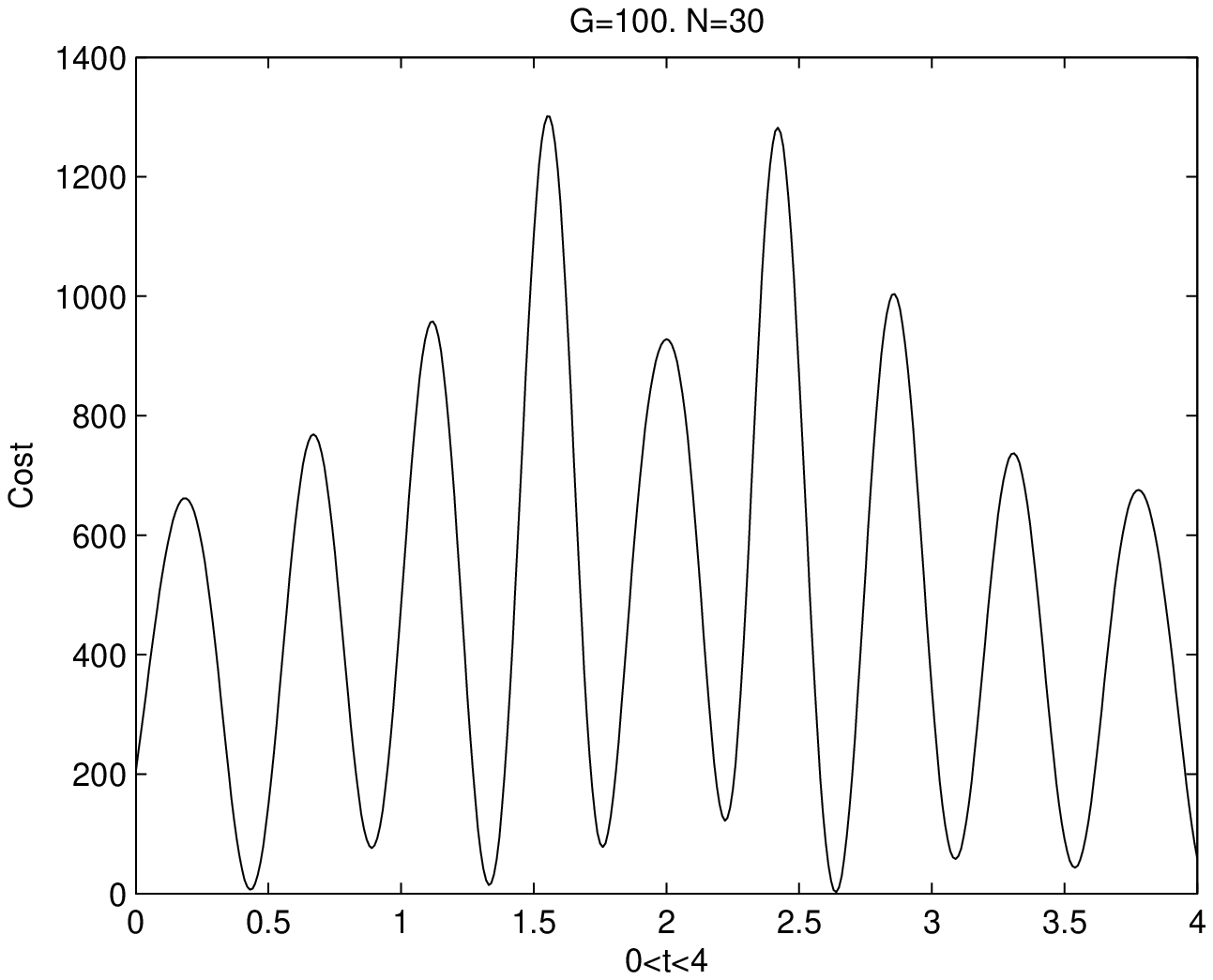}}%
\end{array}
\]

\bigskip

\paragraph{Example 2 : high frequency.}

\bigskip

The initial condition and desired target are specifically as follows. $%
\left( v_{0d},v_{1d}\right) =\left( 0,0\right) $ and with $\left(
k_{o},a_{o},b_{o},\right) =\left( 200,1/2,10000\right) $, for $\left(
x_{1},x_{2}\right) \in \left( 0,1\right) \times \left( 0,1\right) $,
\[
\left\{
\begin{array}{ll}
v_{0}\left( x_{1},x_{2}\right) & =\sum\limits_{j=1}^{G}\left(
\int_{0}^{1}\int_{0}^{1}g_{0}\left( x_{1},x_{2}\right) e_{j}\left(
x_{1},x_{2}\right) dx_{1}dx_{2}\right) e_{j}\left( x_{1},x_{2}\right) \text{
,} \\
v_{1}\left( x_{1},x_{2}\right) & =\sum\limits_{j=1}^{G}\left(
\int_{0}^{1}\int_{0}^{1}g_{1}\left( x_{1},x_{2}\right) e_{j}\left(
x_{1},x_{2}\right) dx_{1}dx_{2}\right) e_{j}\left( x_{1},x_{2}\right) \text{
,} \\
g_{0}\left( x_{1},x_{2}\right) & =e^{-\frac{k_{o}a_{o}}{2}\left(
x_{1}-x_{o1}\right) ^{2}}e^{-\frac{k_{o}b_{o}}{2}\left( x_{2}-x_{o2}\right)
^{2}}\cos \left( k_{o}\left( x_{2}-x_{o2}\right) /2\right) \text{ ,} \\
g_{1}\left( x_{1},x_{2}\right) & =e^{-\frac{k_{o}a_{o}}{2}\left(
x_{1}-x_{o1}\right) ^{2}}e^{-\frac{k_{o}b_{o}}{2}\left( x_{2}-x_{o2}\right)
^{2}} \\
& \quad \left[ k_{o}b_{o}\left( x_{2}-x_{o2}\right) \cos \left( k_{o}\left(
x_{2}-x_{o2}\right) /2\right) \right. \\
& \quad +\left( k_{o}/2+a_{o}\right) \sin \left( k_{o}\left(
x_{2}-x_{o2}\right) /2\right) \\
& \quad \left. -k_{o}a_{o}^{2}\left( x_{1}-x_{o1}\right) ^{2}\sin \left(
k_{o}\left( x_{2}-x_{o2}\right) /2\right) \right] \text{ .}%
\end{array}%
\right.
\]

\noindent Notice that we have chosen as initial data the $G$-first
projections on the basis $\left\{ e_{j}\right\} _{j\geq1}$ of a gaussian
beam $g\left( x_{1},x_{2},t\right) $ such that $g\left( \cdot,t=0\right)
=g_{0}$, $\partial_{t}g\left( \cdot,t=0\right) =g_{1}$\ and which propagates
on the direction $\left( 0,1\right) $.\bigskip

\noindent We take the number of Galerkin mode $G=1000$ and the number of
iterations in the time reversal construction $N=100$.

\bigskip

Below, we plot the graph of the energy of the controlled solution and the
cost of the control function.

\[
\begin{array}{cc}
\raisebox{-0cm}{\includegraphics[ height=4.0918cm, width=5.3065cm
]{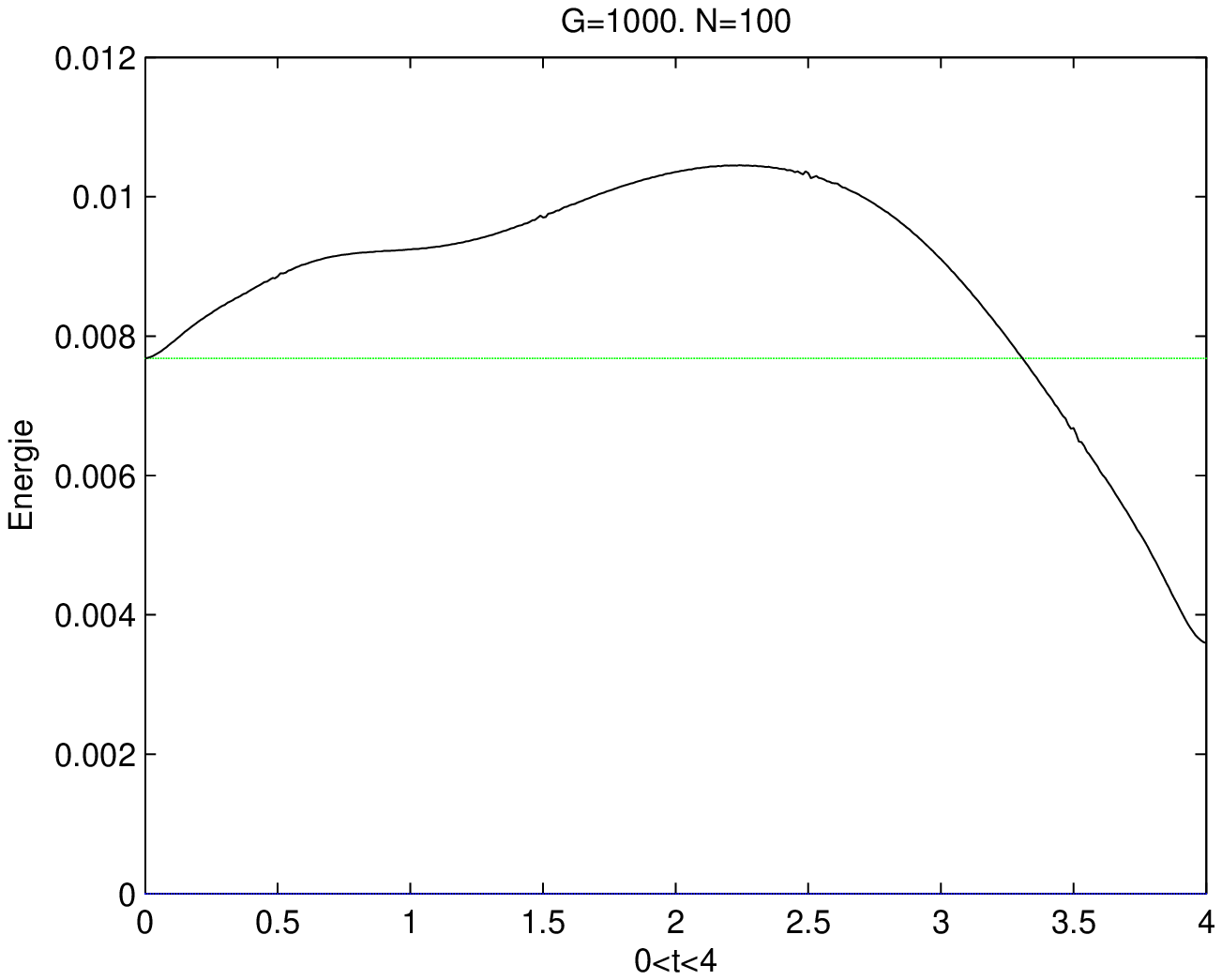}} &
\raisebox{-0cm}{\includegraphics[ height=4.0918cm, width=5.3065cm
]{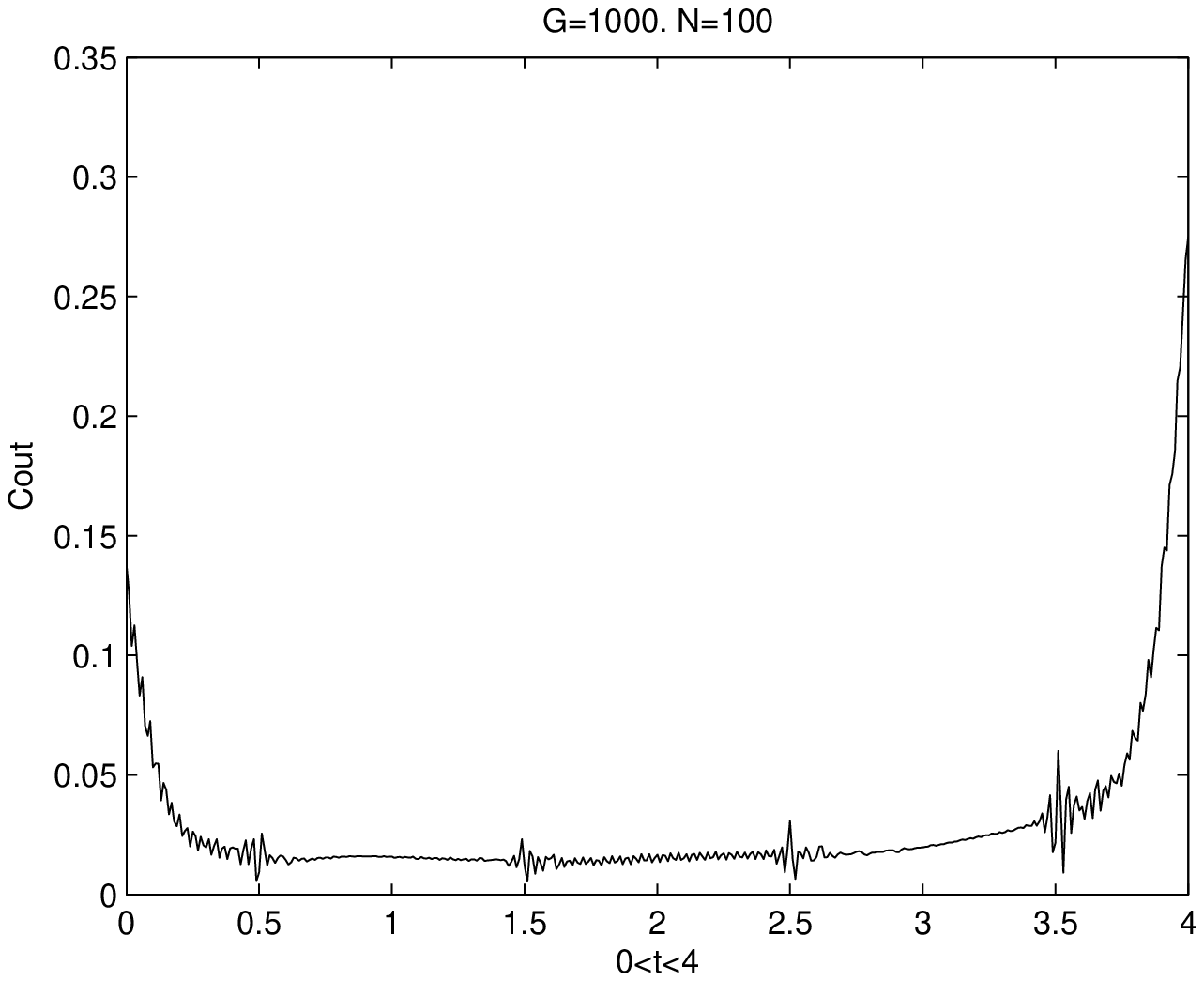}}%
\end{array}
\]

\bigskip

\section{Conclusion}

\bigskip

In this talk, we have considered the wave equation in a bounded domain
(eventually convex). Two kinds of inequality are described when occurs
trapped rays. Applications to control theory are given. First, we link such
kind of estimate with the damped wave equation and its decay rate. Next, we
describe the design of an approximate control function by an iterative time
reversal method. We also provide a numerical simulation in a square domain.
I'm grateful to Prof. Jean-Pierre Puel, the "French-Chinese Summer Institute
on Applied Mathematics" and Fudan University for the kind invitation and the
support of my visit.

\bigskip

\section{Appendix}

\bigskip

In this appendix, we recall most of the material from the works of I.
Kukavica \cite{Ku2} and L. Escauriaza \cite{E} for the elliptic equation and
from the works of G. Lebeau and L. Robbiano \cite{LR} for the wave equation.

\bigskip

\noindent In the original paper dealing with doubling property and frequency
function, N. Garofalo and F.H. Lin \cite{GaL}\ study the monotonicity
property of the following quantity
\[
\frac{r\int_{B_{0,r}}\left| \nabla v\left( y\right) \right| ^{2}dy}{%
\int_{\partial B_{0,r}}\left| v\left( y\right) \right| ^{2}d\sigma\left(
y\right) }\text{ .}
\]

\noindent However, it seems more natural in our context to consider the
monotonicity properties of the frequency function (see \cite{Ze}) defined by
\[
\frac{\int_{B_{0,r}}\left\vert \nabla v\left( y\right) \right\vert
^{2}\left( r^{2}-\left\vert y\right\vert ^{2}\right) dy}{\int_{B_{0,r}}\left%
\vert v\left( y\right) \right\vert ^{2}dy}\text{ .}
\]

\bigskip

\subsection{Monotonicity formula}

\bigskip

\bigskip

Following the ideas of I. Kukavica (\cite{Ku2}, \cite{Ku}, \cite{KN}, see
also \cite{E}, \cite{AE}), one obtains the following three lemmas. Detailed
proofs are given in \cite{Ph3}.

\bigskip

\textbf{Lemma A .-}\qquad \textit{Let }$D\subset \mathbb{R}^{N+1}$, $N\geq 1$%
\textit{, be a connected bounded open set such that }$\overline{%
B_{y_{o},R_{o}}}\subset D$\textit{\ with }$y_{o}\in D$\textit{\ and }$%
R_{o}>0 $\textit{. If }$v=v\left( y\right) \in H^{2}\left( D\right) $\textit{%
\ is a solution of }$\Delta _{y}v=0$\textit{\ in }$D$\textit{, then}
\[
\Phi \left( r\right) =\frac{\int_{B_{y_{o},r}}\left\vert \nabla v\left(
y\right) \right\vert ^{2}\left( r^{2}-\left\vert y-y_{o}\right\vert
^{2}\right) dy}{\int_{B_{y_{o},r}}\left\vert v\left( y\right) \right\vert
^{2}dy}
\]

\noindent \textit{is non-decreasing on }$~0<r<R_{o}$ , \textit{and }%
\[
\frac{d}{dr}\text{ln}\int_{B_{y_{o},r}}\left\vert v\left( y\right)
\right\vert ^{2}dy=\frac{1}{r}\left( N+1+\Phi \left( r\right) \right) \text{
.}
\]

\bigskip

\textbf{Lemma B .-}\qquad\textit{Let }$D\subset\mathbb{R}^{N+1}$, $N\geq 1$%
\textit{, be a connected bounded open set such that }$\overline {%
B_{y_{o},R_{o}}}\subset D$\textit{\ with }$y_{o}\in D$\textit{\ and }$%
R_{o}>0 $\textit{. Let }$r_{1}$\textit{, }$r_{2}$\textit{, }$r_{3}$\textit{\
be three real numbers such that }$0<r_{1}<r_{2}<r_{3}<R_{o}$\textit{. If }$%
v=v\left( y\right) \in H^{2}\left( D\right) $\textit{\ is a solution of }$%
\Delta_{y}v=0$\textit{\ in }$D$\textit{, then}
\[
\int_{B_{y_{o},r_{2}}}\left| v\left( y\right) \right| ^{2}dy\leq\left(
\int_{B_{y_{o},r_{1}}}\left| v\left( y\right) \right| ^{2}dy\right)
^{\alpha}\left( \int_{B_{y_{o},r_{3}}}\left| v\left( y\right) \right|
^{2}dy\right) ^{1-\alpha}\text{ ,}
\]

\noindent\textit{where }$\alpha=\frac{1}{\text{ln}\frac{r_{2}}{r_{1}}}\left(
\frac{1}{\text{ln}\frac{r_{2}}{r_{1}}}+\frac{1}{\text{ln}\frac{r_{3}}{r_{2}}}%
\right) ^{-1}\in\left( 0,1\right) $\textit{.}

\bigskip

\noindent The above two results are still available when we are closed to a
part $\Gamma$ of the boundary $\partial\Omega$ under the homogeneous
Dirichlet boundary condition on $\Gamma$, as follows.

\bigskip

\textbf{Lemma C .-}\qquad\textit{Let }$D\subset\mathbb{R}^{N+1}$, $N\geq 1$%
\textit{, be a connected bounded open set with boundary }$\partial D$\textit{%
. Let }$\Gamma$\textit{\ be a non-empty Lipschitz open subset of }$\partial
D $\textit{. Let }$r_{o}$\textit{, }$r_{1}$\textit{, }$r_{2}$\textit{, }$%
r_{3}$\textit{, }$R_{o}$\textit{\ be five real numbers such that }$%
0<r_{1}<r_{o}<r_{2}<r_{3}<R_{o}$\textit{. Suppose that }$y_{o}\in D$\textit{%
\ satisfies the following three conditions:}

\textit{\quad i). }$B_{y_{o},r}\cap D$\textit{\ is star-shaped with respect
to }$y_{o}\quad\forall r\in\left( 0,R_{o}\right) $\textit{\ ,}

\textit{\quad ii). }$B_{y_{o},r}\subset D\quad\forall r\in\left(
0,r_{o}\right) $\textit{\ ,}

\textit{\quad iii). }$B_{y_{o},r}\cap\partial D\subset\Gamma\quad\forall r\in%
\left[ r_{o},R_{o}\right) $\textit{\ .}

\noindent \textit{If }$v=v\left( y\right) \in H^{2}\left( D\right) $\textit{%
\ is a solution of }$\Delta _{y}v=0$\textit{\ in }$D$\textit{\ and }$v=0$%
\textit{\ on }$\Gamma $,\textit{\ then}
\[
\int_{B_{y_{o},r_{2}}\cap D}\left\vert v\left( y\right) \right\vert
^{2}dy\leq \left( \int_{B_{y_{o},r_{1}}}\left\vert v\left( y\right)
\right\vert ^{2}dy\right) ^{\alpha }\left( \int_{B_{y_{o},r_{3}}\cap
D}\left\vert v\left( y\right) \right\vert ^{2}dy\right) ^{1-\alpha }
\]

\noindent\textit{where }$\alpha=\frac{1}{\text{ln}\frac{r_{2}}{r_{1}}}\left(
\frac{1}{\text{ln}\frac{r_{2}}{r_{1}}}+\frac{1}{\text{ln}\frac{r_{3}}{r_{2}}}%
\right) ^{-1}\in\left( 0,1\right) $\textit{.}

\bigskip

\subsubsection{Proof of Lemma B}

\bigskip

\noindent Let
\[
H\left( r\right) =\int_{B_{y_{o},r}}\left| v\left( y\right) \right| ^{2}dy%
\text{ .}
\]

\noindent By applying Lemma A, we know that
\[
\frac{d}{dr}\text{ln}H\left( r\right) =\frac{1}{r}\left( N+1+\Phi\left(
r\right) \right) \text{ .}
\]

\noindent Next, from the monotonicity property of $\Phi $, one deduces the
following two inequalities

\[
\begin{array}{ll}
\text{ln}\left( \frac{H\left( r_{2}\right) }{H\left( r_{1}\right) }\right) &
=\int_{r_{1}}^{r_{2}}\frac{N+1+\Phi\left( r\right) }{r}dr \\
& \leq\left( N+1+\Phi\left( r_{2}\right) \right) \text{ln}\frac{r_{2}}{r_{1}}%
\text{ ,}%
\end{array}%
\]

\[
\begin{array}{ll}
\text{ln}\left( \frac{H\left( r_{3}\right) }{H\left( r_{2}\right) }\right) &
=\int_{r_{2}}^{r_{3}}\frac{N+1+\Phi\left( r\right) }{r}dr \\
& \geq\left( N+1+\Phi\left( r_{2}\right) \right) \text{ln}\frac{r_{3}}{r_{2}}%
\text{ .}%
\end{array}%
\]

\noindent Consequently,
\[
\frac{\text{ln}\left( \frac{H\left( r_{2}\right) }{H\left( r_{1}\right) }%
\right) }{\text{ln}\frac{r_{2}}{r_{1}}}\leq\left( N+1\right) +\Phi\left(
r_{2}\right) \leq\frac{\text{ln}\left( \frac{H\left( r_{3}\right) }{H\left(
r_{2}\right) }\right) }{\text{ln}\frac{r_{3}}{r_{2}}}\text{ ,}
\]

\noindent and therefore the desired estimate holds
\[
H\left( r_{2}\right) \leq\left( H\left( r_{1}\right) \right) ^{\alpha
}\left( H\left( r_{3}\right) \right) ^{1-\alpha}\text{ ,}
\]

\noindent where $\alpha=\frac{1}{\text{ln}\frac{r_{2}}{r_{1}}}\left( \frac {1%
}{\text{ln}\frac{r_{2}}{r_{1}}}+\frac{1}{\text{ln}\frac{r_{3}}{r_{2}}}%
\right) ^{-1}$.

\bigskip

\subsubsection{Proof of Lemma A}

\bigskip

\noindent We introduce the following two functions $H$ and $D$ for $%
0<r<R_{o} $ :
\[
\begin{array}{ll}
H\left( r\right) = & \int_{B_{y_{o},r}}\left| v\left( y\right) \right| ^{2}dy%
\text{ ,} \\
D\left( r\right) = & \int_{B_{y_{o},r}}\left| \nabla v\left( y\right)
\right| ^{2}\left( r^{2}-\left| y-y_{o}\right| ^{2}\right) dy\text{ .}%
\end{array}%
\]

\bigskip

\noindent First, the derivative of $H\left( r\right)
=\int_{0}^{r}\int_{S^{N}}\left\vert v\left( \rho s+y_{o}\right) \right\vert
^{2}\rho ^{N}d\rho d\sigma \left( s\right) $ is given by $H^{\prime }\left(
r\right) =\int_{\partial B_{y_{o},r}}\left\vert v\left( y\right) \right\vert
^{2}d\sigma \left( y\right) $. Next, recall the Green formula%
\[
\begin{array}{cc}
& \int_{\partial B_{y_{o},r}}\left\vert v\right\vert ^{2}\partial _{\nu
}Gd\sigma \left( y\right) -\int_{\partial B_{y_{o},r}}\partial _{\nu }\left(
\left\vert v\right\vert ^{2}\right) Gd\sigma \left( y\right) \\
& =\int_{B_{y_{o},r}}\left\vert v\right\vert ^{2}\Delta
Gdy-\int_{B_{y_{o},r}}\Delta \left( \left\vert v\right\vert ^{2}\right) Gdy%
\text{ .}%
\end{array}%
\]

\noindent We apply it with $G\left( y\right) =r^{2}-\left| y-y_{o}\right|
^{2}$ where $G_{\left| \partial B_{y_{o},r}\right. }=0$, $\partial_{\nu
}G_{\left| \partial B_{y_{o},r}\right. }=-2r$, and $\Delta G=-2\left(
N+1\right) $. It gives%
\[
\begin{array}{ll}
H^{\prime}\left( r\right) & =\frac{1}{r}\int_{B_{y_{o},r}}\left( N+1\right)
\left| v\right| ^{2}dy+\frac{1}{2r}\int_{B_{y_{o},r}}\Delta\left( \left|
v\right| ^{2}\right) \left( r^{2}-\left| y-y_{o}\right| ^{2}\right) dy \\
& =\frac{N+1}{r}H\left( r\right) +\frac{1}{r}\int_{B_{y_{o},r}}div\left(
v\nabla v\right) \left( r^{2}-\left| y-y_{o}\right| ^{2}\right) dy \\
& =\frac{N+1}{r}H\left( r\right) +\frac{1}{r}\int_{B_{y_{o},r}}\left( \left|
\nabla v\right| ^{2}+v\Delta v\right) \left( r^{2}-\left| y-y_{o}\right|
^{2}\right) dy\text{ .}%
\end{array}%
\]

\noindent Consequently, when $\Delta_{y}v=0$,
\begin{equation}
H^{\prime}\left( r\right) =\frac{N+1}{r}H\left( r\right) +\frac{1}{r}D\left(
r\right) \text{ ,}  \tag{A.1}  \label{A.1}
\end{equation}
\noindent that is $\frac{H^{\prime}\left( r\right) }{H\left( r\right) }=%
\frac{N+1}{r}+\frac{1}{r}\frac{D\left( r\right) }{H\left( r\right) }$ the
second equality in Lemma A.

\bigskip

\noindent Now, we compute the derivative of $D\left( r\right) $.
\begin{equation}
\begin{array}{ll}
D^{\prime }\left( r\right) & =\frac{d}{dr}\left(
r^{2}\int_{0}^{r}\int_{S^{N}}\left\vert \left( \nabla v\right) _{\left\vert
\rho s+y_{o}\right. }\right\vert ^{2}\rho ^{N}d\rho d\sigma \left( s\right)
\right) \\
& \quad -\int_{S^{N}}r^{2}\left\vert \left( \nabla v\right) _{\left\vert
rs+y_{o}\right. }\right\vert ^{2}r^{N}d\sigma \left( s\right) \\
& =2r\int_{0}^{r}\int_{S^{N}}\left\vert \left( \nabla v\right) _{\left\vert
\rho s+y_{o}\right. }\right\vert ^{2}\rho ^{N}d\rho d\sigma \left( s\right)
\\
& =2r\int_{B_{y_{o},r}}\left\vert \nabla v\right\vert ^{2}dy\text{ .}%
\end{array}
\tag{A.2}  \label{A.2}
\end{equation}

\noindent On the other hand, we have by integrations by parts that
\begin{equation}
\begin{array}{ll}
2r\int_{B_{y_{o},r}}\left\vert \nabla v\right\vert ^{2}dy & =\frac{N+1}{r}%
D\left( r\right) +\frac{4}{r}\int_{B_{y_{o},r}}\left\vert \left(
y-y_{o}\right) \cdot \nabla v\right\vert ^{2}dy \\
& \quad -\frac{1}{r}\int_{B_{y_{o},r}}\nabla v\cdot \left( y-y_{o}\right)
\Delta v\left( r^{2}-\left\vert y-y_{o}\right\vert ^{2}\right) dy\text{ .}%
\end{array}
\tag{A.3}  \label{A.3}
\end{equation}

\noindent Therefore,
\[
\begin{array}{ll}
& \quad \left( N+1\right) \int_{B_{y_{o},r}}\left\vert \nabla v\right\vert
^{2}\left( r^{2}-\left\vert y-y_{o}\right\vert ^{2}\right) dy \\
& =2r^{2}\int_{B_{y_{o},r}}\left\vert \nabla v\right\vert
^{2}dy-4\int_{B_{y_{o},r}}\left\vert \left( y-y_{o}\right) \cdot \nabla
v\right\vert ^{2}dy \\
& \quad +2\int_{B_{y_{o},r}}\left( y-y_{o}\right) \cdot \nabla v\Delta
v\left( r^{2}-\left\vert y-y_{o}\right\vert ^{2}\right) dy\text{ ,}%
\end{array}%
\]

\noindent and this is the desired estimate (\ref{A.3}).

\bigskip

\noindent Consequently, from (\ref{A.2}) and (\ref{A.3}), we obtain, when $%
\Delta_{y}v=0$, the following formula%
\begin{equation}
D^{\prime}\left( r\right) =\frac{N+1}{r}D\left( r\right) +\frac{4}{r}%
\int_{B_{y_{o},r}}\left| \left( y-y_{o}\right) \cdot\nabla v\right| ^{2}dy%
\text{ .}  \tag{A.4}  \label{A.4}
\end{equation}

\bigskip

\noindent The computation of the derivative of $\Phi\left( r\right) =\frac{%
D\left( r\right) }{H\left( r\right) }$ gives
\[
\Phi^{\prime}\left( r\right) =\frac{1}{H^{2}\left( r\right) }\left[
D^{\prime}\left( r\right) H\left( r\right) -D\left( r\right) H^{\prime
}\left( r\right) \right] \text{ ,}
\]

\noindent which implies using (\ref{A.1}) and (\ref{A.4}) that
\[
H^{2}\left( r\right) \Phi ^{\prime }\left( r\right) =\frac{1}{r}\left(
4\int_{B_{y_{o},r}}\left\vert \left( y-y_{o}\right) \cdot \nabla
v\right\vert ^{2}dyH\left( r\right) -D^{2}\left( r\right) \right) \geq 0%
\text{ ,}
\]

\noindent indeed, thanks to an integration by parts and using Cauchy-Schwarz
inequality, we have
\[
\begin{array}{ll}
D^{2}\left( r\right) & =4\left( \int_{B_{y_{o},r}}v\nabla v\cdot\left(
y-y_{o}\right) dy\right) ^{2} \\
& \leq4\left( \int_{B_{y_{o},r}}\left| \left( y-y_{o}\right) \cdot\nabla
v\right| ^{2}dy\right) \left( \int_{B_{y_{o},r}}\left| v\right| ^{2}dy\right)
\\
& \leq4\left( \int_{B_{y_{o},r}}\left| \left( y-y_{o}\right) \cdot\nabla
v\right| ^{2}dy\right) H\left( r\right) \text{ .}%
\end{array}%
\]

\noindent Therefore, we have proved the desired monotonicity for $\Phi$ and
this completes the proof of Lemma A.

\bigskip

\subsubsection{Proof of Lemma C}

\bigskip

\noindent Under the assumption $B_{y_{o},r}\cap\partial D\subset\Gamma$ for
any $r\in\left[ r_{o},R_{o}\right) $, we extend $v$ by zero in $\overline{%
B_{y_{o},R_{o}}\left\backslash D\right. }$ and denote by $\overline{v}$ its
extension. Since $v=0$ on $\Gamma$, we have
\[
\left\{
\begin{array}{ll}
\overline{v}=v1_{D} & \text{in}~\overline{B_{y_{o},R_{o}}}\text{ ,} \\
\overline{v}=0 & \text{on}~B_{y_{o},R_{o}}\cap\partial D\text{ ,} \\
\nabla\overline{v}=\nabla v1_{D} & \text{in}~B_{y_{o},R_{o}}\text{ .}%
\end{array}
\right.
\]

\bigskip

\noindent Now, we denote $\Omega_{r}=B_{y_{o},r}\cap D$, when $0<r<R_{o}$.
In particular, $\Omega_{r}=B_{y_{o},r}$, when $0<r<r_{o}$. We introduce the
following three functions:
\[
\begin{array}{ll}
H\left( r\right) = & \int_{\Omega_{r}}\left| v\left( y\right) \right| ^{2}dy%
\text{ ,} \\
D\left( r\right) = & \int_{\Omega_{r}}\left| \nabla v\left( y\right) \right|
^{2}\left( r^{2}-\left| y-y_{o}\right| ^{2}\right) dy\text{ ,}%
\end{array}%
\]

\noindent and
\[
\Phi\left( r\right) =\frac{D\left( r\right) }{H\left( r\right) }\geq0\text{ .%
}
\]

\noindent Our goal is to show that $\Phi$ is a non-decreasing function.
Indeed, we will prove that the following equality holds%
\begin{equation}
\frac{d}{dr}\text{ln}H\left( r\right) =\left( N+1\right) \frac{d}{dr}\text{ln%
}r+\frac{1}{r}\Phi\left( r\right) \text{ .}  \tag{C.1}  \label{C.1}
\end{equation}

\noindent Therefore, from the monotonicity of $\Phi$, we will deduce (in a
similar way than in the proof of Lemma A) that
\[
\frac{\text{ln}\left( \frac{H\left( r_{2}\right) }{H\left( r_{1}\right) }%
\right) }{\text{ln}\frac{r_{2}}{r_{1}}}\leq\left( N+1\right) +\Phi\left(
r_{2}\right) \leq\frac{\text{ln}\left( \frac{H\left( r_{3}\right) }{H\left(
r_{2}\right) }\right) }{\text{ln}\frac{r_{3}}{r_{2}}}\text{ ,}
\]

\noindent and this will imply the desired estimate
\[
\int_{\Omega_{r_{2}}}\left| v\left( y\right) \right| ^{2}dy\leq\left(
\int_{B_{y_{o},r_{1}}}\left| v\left( y\right) \right| ^{2}dy\right)
^{\alpha}\left( \int_{\Omega_{r_{3}}}\left| v\left( y\right) \right|
^{2}dy\right) ^{1-\alpha}\text{ ,}
\]

\noindent where $\alpha=\frac{1}{\text{ln}\frac{r_{2}}{r_{1}}}\left( \frac {1%
}{\text{ln}\frac{r_{2}}{r_{1}}}+\frac{1}{\text{ln}\frac{r_{3}}{r_{2}}}%
\right) ^{-1}$.

\bigskip

\noindent First, we compute the derivative of $H\left( r\right)
=\int_{B_{y_{o},r}}\left\vert \overline{v}\left( y\right) \right\vert ^{2}dy$%
.
\begin{equation}
\begin{array}{ll}
H^{\prime }\left( r\right) & =\int_{S^{N}}\left\vert \overline{v}\left(
rs+y_{o}\right) \right\vert ^{2}r^{N}d\sigma \left( s\right) \\
& =\frac{1}{r}\int_{S^{N}}\left\vert \overline{v}\left( rs+y_{o}\right)
\right\vert ^{2}rs\cdot sr^{N}d\sigma \left( s\right) \\
& =\frac{1}{r}\int_{B_{y_{o},r}}div\left( \left\vert \overline{v}\left(
y\right) \right\vert ^{2}\left( y-y_{o}\right) \right) dy \\
& =\frac{1}{r}\int_{B_{y_{o},r}}\left( \left( N+1\right) \left\vert
\overline{v}\left( y\right) \right\vert ^{2}+\nabla \left\vert \overline{v}%
\left( y\right) \right\vert ^{2}\cdot \left( y-y_{o}\right) \right) dy \\
& =\frac{N+1}{r}H\left( r\right) +\frac{2}{r}\int_{\Omega _{r}}v\left(
y\right) \nabla v\left( y\right) \cdot \left( y-y_{o}\right) dy\text{ .}%
\end{array}
\tag{C.2}  \label{C.2}
\end{equation}

\noindent Next, when $\Delta_{y}v=0$ in $D$ and $v_{\left| \Gamma\right. }=0$%
, we remark that
\begin{equation}
D\left( r\right) =2\int_{\Omega_{r}}v\left( y\right) \nabla v\left( y\right)
\cdot\left( y-y_{o}\right) dy\text{ ,}  \tag{C.3}  \label{C.3}
\end{equation}

\noindent indeed,
\[
\begin{array}{ll}
& \quad \int_{\Omega _{r}}\left\vert \nabla v\right\vert ^{2}\left(
r^{2}-\left\vert y-y_{o}\right\vert ^{2}\right) dy \\
& =\int_{\Omega _{r}}div\left[ v\nabla v\left( r^{2}-\left\vert
y-y_{o}\right\vert ^{2}\right) \right] dy-\int_{\Omega _{r}}vdiv\left[
\nabla v\left( r^{2}-\left\vert y-y_{o}\right\vert ^{2}\right) \right] dy \\
& =-\int_{\Omega _{r}}v\Delta v\left( r^{2}-\left\vert y-y_{o}\right\vert
^{2}\right) dy-\int_{\Omega _{r}}v\nabla v\cdot \nabla \left(
r^{2}-\left\vert y-y_{o}\right\vert ^{2}\right) dy \\
& \quad \text{because on }\partial B_{y_{o},r}\text{, }r=\left\vert
y-y_{o}\right\vert \text{ and }v_{\left\vert \Gamma \right. }=0 \\
& =2\int_{\Omega _{r}}v\nabla v\cdot \left( y-y_{o}\right) dy\quad \text{%
because }\Delta _{y}v=0\text{ in }D\text{ .}%
\end{array}%
\]

\noindent Consequently, from (\ref{C.2}) and (\ref{C.3}), we obtain
\begin{equation}
H^{\prime}\left( r\right) =\frac{N+1}{r}H\left( r\right) +\frac{1}{r}D\left(
r\right) \text{ ,}  \tag{C.4}  \label{C.4}
\end{equation}
\noindent and this is (\ref{C.1}).

\bigskip

\noindent On another hand, the derivative of $D\left( r\right) $ is
\begin{equation}
\begin{array}{ll}
D^{\prime }\left( r\right) & =2r\int_{0}^{r}\int_{S^{N}}\left\vert \left(
\nabla \overline{v}\right) _{\left\vert \rho s+y_{o}\right. }\right\vert
^{2}\rho ^{N}d\rho d\sigma \left( s\right) \\
& =2r\int_{\Omega _{r}}\left\vert \nabla v\left( y\right) \right\vert ^{2}dy%
\text{ .}%
\end{array}
\tag{C.5}  \label{C.5}
\end{equation}

\noindent Here, when $\Delta _{y}v=0$ in $D$ and $v_{\left\vert \Gamma
\right. }=0$, we will remark that
\begin{equation}
\begin{array}{ll}
2r\int_{\Omega _{r}}\left\vert \nabla v\left( y\right) \right\vert ^{2}dy & =%
\frac{N+1}{r}D\left( r\right) +\frac{4}{r}\int_{B_{y_{o},r}}\left\vert
\left( y-y_{o}\right) \cdot \nabla v\left( y\right) \right\vert ^{2}dy \\
& +\frac{1}{r}\int_{\Gamma \cap B_{y_{o},r}}\left\vert \partial _{\nu
}v\right\vert ^{2}\left( r^{2}-\left\vert y-y_{o}\right\vert ^{2}\right)
\left( y-y_{o}\right) \cdot \nu d\sigma \left( y\right)%
\end{array}
\tag{C.6}  \label{C.6}
\end{equation}

\noindent indeed,
\[
\begin{array}{ll}
& \quad \left( N+1\right) \int_{\Omega _{r}}\left\vert \nabla v\right\vert
^{2}\left( r^{2}-\left\vert y-y_{o}\right\vert ^{2}\right) dy \\
& =\int_{\Omega _{r}}div\left( \left\vert \nabla v\right\vert ^{2}\left(
r^{2}-\left\vert y-y_{o}\right\vert ^{2}\right) \left( y-y_{o}\right)
\right) dy \\
& \quad -\int_{\Omega _{r}}\nabla \left( \left\vert \nabla v\right\vert
^{2}\left( r^{2}-\left\vert y-y_{o}\right\vert ^{2}\right) \right) \cdot
\left( y-y_{o}\right) dy \\
& =\int_{\Gamma \cap B_{y_{o},r}}\left\vert \nabla v\right\vert ^{2}\left(
r^{2}-\left\vert y-y_{o}\right\vert ^{2}\right) \left( y-y_{o}\right) \cdot
\nu d\sigma \left( y\right) \\
& \quad -\int_{\Omega _{r}}\partial _{y_{i}}\left( \left\vert \nabla
v\right\vert ^{2}\left( r^{2}-\left\vert y-y_{o}\right\vert ^{2}\right)
\right) \left( y_{i}-y_{oi}\right) dy \\
& =\int_{\Gamma \cap B_{y_{o},r}}\left\vert \nabla v\right\vert ^{2}\left(
r^{2}-\left\vert y-y_{o}\right\vert ^{2}\right) \left( y-y_{o}\right) \cdot
\nu d\sigma \left( y\right) \\
& \quad -\int_{\Omega _{r}}2\nabla v\partial _{y_{i}}\nabla v\left(
r^{2}-\left\vert y-y_{o}\right\vert ^{2}\right) \left( y_{i}-y_{oi}\right) dy
\\
& \quad +2\int_{\Omega _{r}}\left\vert \nabla v\right\vert ^{2}\left\vert
y-y_{o}\right\vert ^{2}dy\text{ ,}%
\end{array}%
\]

\[
\begin{array}{ll}
\text{ and \qquad} & \quad-\int_{\Omega_{r}}\partial_{y_{j}}v\partial
_{y_{i}y_{j}}^{2}v\left( r^{2}-\left| y-y_{o}\right| ^{2}\right) \left(
y_{i}-y_{oi}\right) dy \\
& =-\int_{\Omega_{r}}\partial_{y_{j}}\left( \left( y_{i}-y_{oi}\right)
\partial_{y_{j}}v\partial_{y_{i}}v\left( r^{2}-\left| y-y_{o}\right|
^{2}\right) \right) dy \\
& \quad+\int_{\Omega_{r}}\partial_{y_{j}}\left( y_{i}-y_{oi}\right)
\partial_{y_{j}}v\partial_{y_{i}}v\left( r^{2}-\left| y-y_{o}\right|
^{2}\right) dy \\
& \quad+\int_{\Omega_{r}}\left( y_{i}-y_{oi}\right)
\partial_{y_{j}}^{2}v\partial_{y_{i}}v\left( r^{2}-\left| y-y_{o}\right|
^{2}\right) dy \\
& \quad+\int_{\Omega_{r}}\left( y_{i}-y_{oi}\right)
\partial_{y_{j}}v\partial_{y_{i}}v\partial_{y_{j}}\left( r^{2}-\left|
y-y_{o}\right| ^{2}\right) dy \\
& =-\int_{\Gamma\cap B_{y_{o},r}}\nu_{_{j}}\left( \left( y_{i}-y_{oi}\right)
\partial_{y_{j}}v\partial_{y_{i}}v\left( r^{2}-\left| y-y_{o}\right|
^{2}\right) \right) d\sigma\left( y\right) \\
& \quad+\int_{\Omega_{r}}\left| \nabla v\right| ^{2}\left( r^{2}-\left|
y-y_{o}\right| ^{2}\right) dy \\
& \quad+0\quad\text{because }\Delta_{y}v=0\text{ in }D\text{ } \\
& \quad-\int_{\Omega_{r}}2\left| \left( y-y_{o}\right) \cdot\nabla v\right|
^{2}dy\text{ .}%
\end{array}%
\]

\noindent Therefore, when $\Delta _{y}v=0$ in $D$, we have%
\[
\begin{array}{ll}
& \quad \left( N+1\right) \int_{\Omega _{r}}\left\vert \nabla v\right\vert
^{2}\left( r^{2}-\left\vert y-y_{o}\right\vert ^{2}\right) dy \\
& =\int_{\Gamma \cap B_{y_{o},r}}\left\vert \nabla v\right\vert ^{2}\left(
r^{2}-\left\vert y-y_{o}\right\vert ^{2}\right) \left( y-y_{o}\right) \cdot
\nu d\sigma \left( y\right) \\
& \quad -2\int_{\Gamma \cap B_{y_{o},r}}\partial _{y_{j}}v\nu _{_{j}}\left(
\left( y_{i}-y_{oi}\right) \partial _{y_{i}}v\right) \left( r^{2}-\left\vert
y-y_{o}\right\vert ^{2}\right) d\sigma \left( y\right) \\
& \quad +2r^{2}\int_{\Omega _{r}}\left\vert \nabla u\right\vert
^{2}dy-4\int_{\Omega _{r}}\left\vert \left( y-y_{o}\right) \cdot \nabla
v\right\vert ^{2}dy\text{ .}%
\end{array}%
\]

\noindent By using the fact that $v_{\left\vert \Gamma \right. }=0$, we get $%
\nabla v=\left( \nabla v\cdot \nu \right) \nu $ on $\Gamma $ and we deduce
that
\[
\begin{array}{ll}
& \quad \left( N+1\right) \int_{\Omega _{r}}\left\vert \nabla v\right\vert
^{2}\left( r^{2}-\left\vert y-y_{o}\right\vert ^{2}\right) dy \\
& =-\int_{\Gamma \cap B_{y_{o},r}}\left\vert \partial _{\nu }v\right\vert
^{2}\left( r^{2}-\left\vert y-y_{o}\right\vert ^{2}\right) \left(
y-y_{o}\right) \cdot \nu d\sigma \left( y\right) \\
& \quad +2r^{2}\int_{\Omega _{r}}\left\vert \nabla v\right\vert
^{2}dy-4\int_{\Omega _{r}}\left\vert \left( y-y_{o}\right) \cdot \nabla
v\right\vert ^{2}dy\text{ ,}%
\end{array}%
\]

\noindent and this is (\ref{C.6}).

\bigskip

\noindent Consequently, from (\ref{C.5}) and (\ref{C.6}), when $%
\Delta_{y}v=0 $ in $D$ and $v_{\left| \Gamma\right. }=0$, we have%
\begin{equation}
\begin{array}{ll}
D^{\prime}\left( r\right) & =\frac{N+1}{r}D\left( r\right) +\frac{4}{r}%
\int_{\Omega_{r}}\left| \left( y-y_{o}\right) \cdot\nabla v\left( y\right)
\right| ^{2}dy \\
& \quad+\frac{1}{r}\int_{\Gamma\cap B_{y_{o},r}}\left| \partial_{\nu
}v\right| ^{2}\left( r^{2}-\left| y-y_{o}\right| ^{2}\right) \left(
y-y_{o}\right) \cdot\nu d\sigma\left( y\right) \text{ .}%
\end{array}
\tag{C.7}  \label{C.7}
\end{equation}

\bigskip

\noindent The computation of the derivative of $\Phi\left( r\right) =\frac{%
D\left( r\right) }{H\left( r\right) }$ gives
\[
\Phi^{\prime}\left( r\right) =\frac{1}{H^{2}\left( r\right) }\left[
D^{\prime}\left( r\right) H\left( r\right) -D\left( r\right) H^{\prime
}\left( r\right) \right] \text{ ,}
\]

\noindent which implies from (\ref{C.4}) and (\ref{C.7}), that%
\[
\begin{array}{ll}
H^{2}\left( r\right) \Phi ^{\prime }\left( r\right) & =\frac{1}{r}\left(
4\int_{\Omega _{r}}\left\vert \left( y-y_{o}\right) \cdot \nabla v\left(
y\right) \right\vert ^{2}dy~H\left( r\right) -D^{2}\left( r\right) \right)
\\
& \quad +\frac{H\left( r\right) }{r}\int_{\Gamma \cap B_{y_{o},r}}\left\vert
\partial _{\nu }v\right\vert ^{2}\left( r^{2}-\left\vert y-y_{o}\right\vert
^{2}\right) \left( y-y_{o}\right) \cdot \nu d\sigma \left( y\right)%
\end{array}%
\]

\noindent Thanks to (\ref{C.3}) and Cauchy-Schwarz inequality, we obtain
that
\[
0\leq4\int_{\Omega_{r}}\left| \left( y-y_{o}\right) \cdot\nabla v\left(
y\right) \right| ^{2}dy~H\left( r\right) -D^{2}\left( r\right) \text{ .}
\]

\noindent The inequality $0\leq \left( y-y_{o}\right) \cdot \nu $ on $\Gamma
$ holds when $B_{y_{o},r}\cap D$ is star-shaped with respect to $y_{o}$ for
any $r\in \left( 0,R_{o}\right) $. Therefore, we get the desired
monotonicity for $\Phi $ which completes the proof of Lemma C.

\bigskip

\subsection{Quantitative unique continuation property for the Laplacian}

\bigskip

\noindent Let $D\subset\mathbb{R}^{N+1}$, $N\geq1$, be a connected
bounded open set with boundary $\partial D$. Let $\Gamma$ be a
non-empty Lipschitz open part of $\partial D$. We consider the
Laplacian in $D$, with a homogeneous Dirichlet boundary condition on
$\Gamma\subset\partial\Omega$:
\begin{equation}
\left\{
\begin{array}{rl}
\Delta_{y}v=0 & \quad\text{in}~D\text{ ,} \\
v=0 & \quad\text{on}~\Gamma\text{ ,} \\
v=v\left( y\right) \in H^{2}\left( D\right) & \text{ .}%
\end{array}
\right.  \tag{D.1}  \label{D.1}
\end{equation}

\noindent The goal of this section is to describe interpolation inequalities
associated to solutions $v$ of (\ref{D.1}).

\bigskip

\textbf{Theorem D .-}\qquad\textit{Let }$\omega$\textit{\ be a non-empty
open subset of }$D$\textit{. Then, for any }$D_{1}\subset D$\textit{\ such
that }$\partial D_{1}\cap\partial D\Subset\Gamma$\textit{\ and }$\overline{%
D_{1}}\left\backslash \left( \Gamma\cap\partial D_{1}\right) \right. \subset
D$\textit{, there exist }$C>0$\textit{\ and }$\mu\in\left( 0,1\right) $%
\textit{\ such that for any }$v$\textit{\ solution of }(\ref{D.1})\textit{,
we have}
\[
\int_{D_{1}}\left| v\left( y\right) \right| ^{2}dy\leq C\left(
\int_{\omega}\left| v\left( y\right) \right| ^{2}dy\right) ^{\mu}\left(
\int_{D}\left| v\left( y\right) \right| ^{2}dy\right) ^{1-\mu}\text{ .}
\]

\bigskip

\noindent Or in a equivalent way by a minimization technique,

\bigskip

\textbf{Theorem D' .-}\qquad\textit{Let }$\omega$\textit{\ be a non-empty
open subset of }$D$\textit{. Then, for any }$D_{1}\subset D$\textit{\ such
that }$\partial D_{1}\cap\partial D\Subset\Gamma$\textit{\ and }$\overline{%
D_{1}}\left\backslash \left( \Gamma\cap\partial D_{1}\right) \right. \subset
D$\textit{, there exist }$C>0$\textit{\ and }$\mu\in\left( 0,1\right) $%
\textit{\ such that for any }$v$\textit{\ solution of }(\ref{D.1})\textit{,
we have}
\[
\int_{D_{1}}\left| v\left( y\right) \right| ^{2}dy\leq C\left( \frac {1}{%
\varepsilon}\right) ^{\frac{1-\mu}{\mu}}\int_{\omega}\left| v\left( y\right)
\right| ^{2}dy+\varepsilon\int_{D}\left| v\left( y\right) \right|
^{2}dy\quad\forall\varepsilon>0\text{ .}
\]

\bigskip

Proof of Theorem D .- We divide the proof into two steps.

\bigskip

Step 1 .- We apply Lemma B, and use a standard argument (see e.g., \cite{Ro}%
) which consists to construct a sequence of balls chained along a curve.
More precisely, we claim that for any non-empty compact sets in $D$, $K_{1}$
and $K_{2}$, such that meas$\left( K_{1}\right) >0$, there exists $\mu \in
\left( 0,1\right) $ such that for any $v=v\left( y\right) \in H^{2}\left(
D\right) $, solution of $\Delta _{y}v=0$\ in $D$, we have
\begin{equation}
\int_{K_{2}}\left\vert v\left( y\right) \right\vert ^{2}dy\leq \left(
\int_{K_{1}}\left\vert v\left( y\right) \right\vert ^{2}dy\right) ^{\mu
}\left( \int_{D}\left\vert v\left( y\right) \right\vert ^{2}dy\right)
^{1-\mu }\text{ .}  \tag{D.2}  \label{D.2}
\end{equation}

\bigskip

Step 2 .- We apply Lemma C, and choose $y_{o}$ in a neighborhood of the part
$\Gamma $ such that the conditions \textit{i, ii, iii}, hold. Next, by an
adequate partition of $D$, we deduce from (\ref{D.2}) that for any $%
D_{1}\subset D$ such that $\partial D_{1}\cap \partial D\Subset \Gamma $ and
$\overline{D_{1}}\left\backslash \left( \Gamma \cap \partial D_{1}\right)
\right. \subset D$, there exist $C>0$\ and $\mu \in \left( 0,1\right) $\
such that for any $v=v\left( y\right) \in H^{2}\left( D\right) $\ such that $%
\Delta _{y}v=0$\ on $D$\ and $v=0$ on $\Gamma $, we have
\[
\int_{D_{1}}\left\vert v\left( y\right) \right\vert ^{2}dy\leq C\left(
\int_{\omega }\left\vert v\left( y\right) \right\vert ^{2}dy\right) ^{\mu
}\left( \int_{D}\left\vert v\left( y\right) \right\vert ^{2}dy\right)
^{1-\mu }\text{ .}
\]%
This completes the proof.

\bigskip

\subsection{Quantitative unique continuation property for the elliptic
operator $\partial _{t}^{2}+\Delta $}

\bigskip

In this section, we present the following result.

\bigskip

\textbf{Theorem E .-}\qquad\textit{Let }$\Omega$\textit{\ be a bounded open
set in }$\mathbb{R}^{n}$\textit{, }$n\geq1$\textit{, either convex or }$%
C^{2} $\textit{\ and connected. We choose }$T_{2}>T_{1}$\textit{\ and }$%
\delta \in\left( 0,\left( T_{2}-T_{1}\right) /2\right) $\textit{. Let }$f\in
L^{2}\left( \Omega\times\left( T_{1},T_{2}\right) \right) $. \textit{We
consider the elliptic operator of second order in }$\Omega\times\left(
T_{1},T_{2}\right) $\textit{\ with a homogeneous Dirichlet boundary
condition on }$\partial\Omega\times\left( T_{1},T_{2}\right) $\textit{,}
\begin{equation}
\left\{
\begin{array}{ll}
\partial_{t}^{2}w+\Delta w=f & \text{\textit{in}}~\Omega\times\left(
T_{1},T_{2}\right) \text{ ,} \\
w=0 & \text{\textit{on}}~\partial\Omega\times\left( T_{1},T_{2}\right) \text{
,} \\
w=w\left( x,t\right) \in H^{2}\left( \Omega\times\left( T_{1},T_{2}\right)
\right) & \text{ .}%
\end{array}
\right.  \tag{E.1}  \label{E.1}
\end{equation}

\noindent \textit{Then, for any }$\varphi \in C_{0}^{\infty }\left( \Omega
\times \left( T_{1},T_{2}\right) \right) $\textit{, }$\varphi \neq 0$\textit{%
, there exist }$C>0$\textit{\ and }$\mu \in \left( 0,1\right) $\textit{\
such that for any }$w$\textit{\ solution of (}\ref{E.1}\textit{), we have}
\[
\begin{array}{ll}
& \quad \int_{T_{1}+\delta }^{T_{2}-\delta }\int_{\Omega }\left\vert w\left(
x,t\right) \right\vert ^{2}dxdt \\
& \leq C\left( \int_{T_{1}}^{T_{2}}\int_{\Omega }\left\vert w\left(
x,t\right) \right\vert ^{2}dxdt\right) ^{1-\mu } \\
& \quad \quad \left( \int_{T_{1}}^{T_{2}}\int_{\Omega }\left\vert \varphi
w\left( x,t\right) \right\vert ^{2}dxdt+\int_{T_{1}}^{T_{2}}\int_{\Omega
}\left\vert f\left( x,t\right) \right\vert ^{2}dxdt\right) ^{\mu }\text{ .}%
\end{array}%
\]

\bigskip

Proof .- First, by a difference quotient technique and a standard extension
at $\Omega \times \left\{ T_{1},T_{2}\right\} $, we check the existence of a
solution $u\in H^{2}\left( \Omega \times \left( T_{1},T_{2}\right) \right) $
solving%
\[
\left\{
\begin{array}{ll}
\partial _{t}^{2}u+\Delta u=f & \quad \text{in}~\Omega \times \left(
T_{1},T_{2}\right) \text{ ,} \\
u=0 & \quad \text{on}~\partial \Omega \times \left( T_{1},T_{2}\right) \cup
\Omega \times \left\{ T_{1},T_{2}\right\} \text{ ,}%
\end{array}%
\right.
\]

\noindent such that
\[
\left\| u\right\| _{H^{2}\left( \Omega\times\left( T_{1},T_{2}\right)
\right) }\leq c\left\| f\right\| _{L^{2}\left( \Omega\times\left(
T_{1},T_{2}\right) \right) }\text{ ,}
\]

\noindent for some $c>0$ only depending on $\left( \Omega,T_{1},T_{2}\right)
$. Next, we apply Theorem D with $D=\Omega\times\left( T_{1},T_{2}\right) $,
$\Omega\times\left( T_{1}+\delta,T_{2}-\delta\right) \subset D_{1}$, $%
y=\left( x,t\right) $, $\Delta_{y}=\partial_{t}^{2}+\Delta$, and $v=w-u$.

\bigskip

\subsection{Application to the wave equation}

\bigskip

From the idea of L. Robbiano \cite{Ro2} which consists to use an
interpolation inequality of H\"{o}lder type for the elliptic operator $%
\partial_{t}^{2}+\Delta$ and the Fourier-Bros-Iagolnitzer transform
introduced by G. Lebeau and L. Robbiano \cite{LR}, we obtain the following
estimate of logarithmic type.

\bigskip

\textbf{Theorem F .-}\qquad \textit{Let }$\Omega $\textit{\ be a bounded
open set in }$\mathbb{R}^{n}$\textit{, }$n\geq 1$\textit{, either convex or }%
$C^{2}$\textit{\ and connected. Let }$\omega $\textit{\ be a non-empty open
subset in }$\Omega $\textit{. Then, for any }$\beta \in \left( 0,1\right) $%
\textit{\ and }$k\in \mathbb{N}^{\ast }$\textit{, there exist }$C>0$\textit{%
\ and }$T>0$\textit{\ such that for any solution }$u$\textit{\ of}%
\[
\left\{
\begin{array}{rl}
\partial _{t}^{2}u-\Delta u=0 & \quad \text{\textit{in}~}\Omega \times
\left( 0,T\right) \ \text{,} \\
u=0 & \quad \text{\textit{on}~}\partial \Omega \times \left( 0,T\right)
\text{ ,} \\
\left( u,\partial _{t}u\right) \left( \cdot ,0\right) =\left(
u_{0},u_{1}\right) & \text{ ,}%
\end{array}%
\right.
\]%
\textit{with non-identically zero initial data }$\left( u_{0},u_{1}\right)
\in D\left( A^{k-1}\right) $\textit{, we have}%
\[
\left\Vert \left( u_{0},u_{1}\right) \right\Vert _{D\left( A^{k-1}\right)
}\leq Ce^{\left( C\frac{\left\Vert \left( u_{0},u_{1}\right) \right\Vert
_{D\left( A^{k-1}\right) }}{\left\Vert \left( u_{0},u_{1}\right) \right\Vert
_{L^{2}\left( \Omega \right) \times H^{-1}\left( \Omega \right) }}\right)
^{1/\left( \beta k\right) }}~\left\Vert u\right\Vert _{L^{2}\left( \omega
\times \left( 0,T\right) \right) }\text{ .}
\]
(Here, $D\left( A^{0}\right)=H_{0}^{1}\left( \Omega \right) \times
L^{2}\left( \Omega \right) $).
\bigskip

Proof .- First, recall that with a standard energy method, we have that
\begin{equation}
\forall t\in \mathbb{R}\qquad \left\Vert \left( u_{0},u_{1}\right)
\right\Vert _{H_{0}^{1}\left( \Omega \right) \times L^{2}\left( \Omega
\right) }^{2}=\int_{\Omega }\left( \left\vert \partial _{t}u\left(
x,t\right) \right\vert ^{2}+\left\vert \nabla u\left( x,t\right) \right\vert
^{2}\right) dx\text{ ,}  \tag{F.1}  \label{F.1}
\end{equation}%
and there exists a constant $c>0$ such that for all $T\geq 1$,
\begin{equation}
T\left\Vert \left( u_{0},u_{1}\right) \right\Vert _{L^{2}\left( \Omega
\right) \times H^{-1}\left( \Omega \right) }^{2}\leq
c\int_{0}^{T}\int_{\Omega }\left\vert u\left( x,t\right) \right\vert ^{2}dx%
\text{ .}  \tag{F.2}  \label{F.2}
\end{equation}%
Next, let $\beta \in \left( 0,1\right) $, $k\in \mathbb{N}^{\ast }$, and
choose $N\in \mathbb{N}^{\ast }$ such that $0<\beta +\frac{1}{2N}<1$ and $%
2N>k$. Put $\gamma =1-\frac{1}{2N}$. For any $\lambda \geq 1$, the function $%
F_{\lambda }(z)=\frac{1}{2\pi }\int_{\mathbb{R}}e^{iz\tau }e^{-\left( \frac{%
\tau }{\lambda ^{\gamma }}\right) ^{2N}}d\tau $ is holomorphic in $\mathbb{C}
$, and there exists four positive constants $C_{o}$, $c_{0}$, $c_{1}$ and $%
c_{2}$ (independent on $\lambda $) such that
\begin{equation}
\left\{
\begin{array}{ll}
\forall z\in \mathbb{C}\quad \left\vert F_{\lambda }(z)\right\vert \leq
C_{o}\lambda ^{\gamma }e^{c_{0}\lambda \left\vert \text{Im}z\right\vert
^{1/\gamma }}\text{ ,} &  \\
\left\vert \text{Im}z\right\vert \leq c_{2}\left\vert \text{Re}z\right\vert
\Rightarrow \left\vert F_{\lambda }(z)\right\vert \leq C_{o}\lambda ^{\gamma
}e^{-c_{1}\lambda \left\vert \text{Re}z\right\vert ^{1/\gamma }}\text{ ,} &
\end{array}%
\right.  \tag{F.3}  \label{F.3}
\end{equation}%
(see \cite{LR}).

\bigskip

\noindent Now, let $s,\ell _{o}\in \mathbb{R}$, we introduce the following
Fourier-Bros-Iagolnitzer transformation in \cite{LR}:%
\begin{equation}
W_{\ell _{o},\lambda }(x,s)=\int_{\mathbb{R}}F_{\lambda }(\ell _{o}+is-\ell
)\Phi (\ell )u(x,\ell )d\ell \text{ ,}  \tag{F.4}  \label{F.4}
\end{equation}%
where $\Phi \in C_{0}^{\infty }(\mathbb{R})$. As $u$ is solution of the wave
equation, $W_{\ell _{o},\lambda }$ satisfies:%
\begin{equation}
\left\{
\begin{array}{l}
\partial _{s}^{2}W_{\ell _{o},\lambda }(x,s)+\Delta W_{\ell _{o},\lambda
}(x,s) \\
\quad =\int_{\mathbb{R}}-F_{\lambda }(\ell _{o}+is-\ell )\left[ \Phi
^{\prime \prime }(\ell )u(x,\ell )+2\Phi ^{\prime }(\ell )\partial
_{t}u(x,\ell )\right] d\ell \text{ ,} \\
W_{\ell _{o},\lambda }(x,s)=0\quad \text{for }x\in \partial \Omega \text{ ,}
\\
W_{\ell _{o},\lambda }(x,0)=\left( F_{\lambda }\ast \Phi u(x,\cdot )\right)
(\ell _{o})\quad \text{for }x\in \Omega \text{ .}%
\end{array}%
\right.  \tag{F.5}  \label{F.5}
\end{equation}%
In another hand, we also have for any $T>0$,
\begin{equation}
\begin{array}{ll}
\left\Vert \Phi u\left( x,\cdot \right) \right\Vert _{L^{2}\left( \left(
\frac{T}{2}-1,\frac{T}{2}+1\right) \right) } & \leq \left\Vert \Phi
u(x,\cdot )-F_{\lambda }\ast \Phi u(x,\cdot )\right\Vert _{L^{2}\left(
\left( \frac{T}{2}-1,\frac{T}{2}+1\right) \right) } \\
& \quad +\left\Vert F_{\lambda }\ast \Phi u(x,\cdot )\right\Vert
_{L^{2}\left( \left( \frac{T}{2}-1,\frac{T}{2}+1\right) \right) } \\
& \leq \left\Vert \Phi u(x,\cdot )-F_{\lambda }\ast \Phi u(x,\cdot
)\right\Vert _{L^{2}\left( \mathbb{R}\right) } \\
& \quad +\left( \int_{t\in \left( \frac{T}{2}-1,\frac{T}{2}+1\right)
}\left\vert W_{t,\lambda }(x,0)\right\vert ^{2}dt\right) ^{1/2}\text{ .}%
\end{array}
\tag{F.6}  \label{F.6}
\end{equation}%
Denoting $\mathcal{F}\left( f\right) $ the Fourier transform of $f$, by
using Parseval equality and $\mathcal{F}\left( F_{\lambda }\right) \left(
\tau \right) =e^{-\left( \frac{\tau }{\lambda ^{\gamma }}\right) ^{2N}}$,
one obtain%
\begin{equation}
\begin{array}{ll}
& \quad \left\Vert \Phi u(x,\cdot )-F_{\lambda }\ast \Phi u(x,\cdot
)\right\Vert _{L^{2}\left( \mathbb{R}\right) } \\
& =\frac{1}{\sqrt{2\pi }}\left\Vert \mathcal{F}\left( \Phi u(x,\cdot
)-F_{\lambda }\ast \Phi u(x,\cdot )\right) \right\Vert _{L^{2}\left( \mathbb{%
R}\right) } \\
& =\frac{1}{\sqrt{2\pi }}\left( \int_{\mathbb{R}}\left\vert \left(
1-e^{-\left( \frac{\tau }{\lambda ^{\gamma }}\right) ^{2N}}\right) \mathcal{F%
}\left( \Phi u(x,\cdot )\right) \left( \tau \right) \right\vert ^{2}d\tau
\right) ^{1/2} \\
& \leq C\left( \int_{\mathbb{R}}\left\vert \left( \frac{\tau }{\lambda
^{\gamma }}\right) ^{k}\mathcal{F}\left( \Phi u(x,\cdot )\right) \left( \tau
\right) \right\vert ^{2}d\tau \right) ^{1/2}\quad \text{because }k<2N \\
& \leq C\frac{1}{\lambda ^{\beta k}}\left( \int_{\mathbb{R}}\left\vert
\mathcal{F}\left( \partial _{t}^{k}\left( \Phi u(x,\cdot )\right) \right)
\left( \tau \right) \right\vert ^{2}d\tau \right) ^{1/2}\quad \text{because }%
\beta <\gamma \\
& \leq C\frac{1}{\lambda ^{\beta k}}\left\Vert \partial _{t}^{k}\left( \Phi
u(x,\cdot )\right) \right\Vert _{L^{2}\left( \mathbb{R}\right) }\text{ .}%
\end{array}
\tag{F.7}  \label{F.7}
\end{equation}%
Therefore, from (\ref{F.6}) and (\ref{F.7}), one gets%
\begin{equation}
\begin{array}{ll}
& \quad \left\Vert \Phi u\left( x,\cdot \right) \right\Vert _{L^{2}\left(
\left( \frac{T}{2}-1,\frac{T}{2}+1\right) \right) } \\
& \leq C\frac{1}{\lambda ^{\beta k}}\left\Vert \partial _{t}^{k}\left( \Phi
u(x,\cdot )\right) \right\Vert _{L^{2}\left( \mathbb{R}\right) }+\left(
\int_{t\in \left( \frac{T}{2}-1,\frac{T}{2}+1\right) }\left\vert
W_{t,\lambda }(x,0)\right\vert ^{2}dt\right) ^{1/2}\text{ .}%
\end{array}
\tag{F.8}  \label{F.8}
\end{equation}

\bigskip

\noindent Now, recall that from the Cauchy's theorem we have:

\bigskip

\textbf{Proposition 1 .-}\qquad \textit{Let }$f$\textit{\ be a holomorphic
function in a domain }$D\subset \mathbb{C}$\textit{. Let }$a,b>0$\textit{, }$%
z\in \mathbb{C}$\textit{. We suppose that}
\[
D_{o}=\left\{ \left( x,y\right) \in \mathbb{R}^{2}\simeq \mathbb{C}%
\left\backslash {}\right. \left\vert x-\text{Re}z\right\vert \leq
a,~\left\vert y-\text{Im}z\right\vert \leq b\right\} \subset D\text{ ,}
\]%
\textit{then}
\[
f\left( z\right) =\frac{1}{\pi ab}\int \int_{\left\vert \frac{x-\text{Re}z}{a%
}\right\vert ^{2}+\left\vert \frac{y-\text{Im}z}{b}\right\vert ^{2}\leq
1}f\left( x+iy\right) dxdy\text{ .}
\]

\bigskip

\noindent Choosing $z=t\in \left( \frac{T}{2}-1,\frac{T}{2}+1\right) \subset
\mathbb{R}$ and $x+iy=\ell _{o}+is$, we deduce that
\begin{equation}
\begin{array}{ll}
\left\vert W_{t,\lambda }(x,0)\right\vert & \leq \frac{1}{\pi ab}%
\int_{\left\vert \ell _{o}-t\right\vert \leq a}\int_{\left\vert s\right\vert
\leq b}\left\vert W_{\ell _{o}+is,\lambda }(x,0)\right\vert d\ell _{o}ds \\
& \leq \frac{1}{\pi ab}\int_{\left\vert \ell _{o}-t\right\vert \leq
a}\int_{\left\vert s\right\vert \leq b}\left\vert W_{\ell _{o},\lambda
}(x,s)\right\vert dsd\ell _{o} \\
& \leq \frac{2}{\pi \sqrt{ab}}\left( \int_{\left\vert \ell _{o}-t\right\vert
\leq a}\int_{\left\vert s\right\vert \leq b}\left\vert W_{\ell _{o},\lambda
}(x,s)\right\vert ^{2}dsd\ell _{o}\right) ^{1/2}\text{ ,}%
\end{array}
\tag{F.9}  \label{F.9}
\end{equation}%
and with $a=2b=1$,
\begin{equation}
\begin{array}{ll}
& \quad \int_{t\in \left( \frac{T}{2}-1,\frac{T}{2}+1\right) }\left\vert
W_{t,\lambda }(x,0)\right\vert ^{2}dt \\
& \leq \int_{t\in \left( \frac{T}{2}-1,\frac{T}{2}+1\right) }\left(
\int_{\left\vert \ell _{o}-t\right\vert \leq 1}\int_{\left\vert s\right\vert
\leq 1/2}\left\vert W_{\ell _{o},\lambda }(x,s)\right\vert ^{2}dsd\ell
_{o}\right) dt \\
& \leq \int_{t\in \left( \frac{T}{2}-1,\frac{T}{2}+1\right) }\int_{\ell
_{o}\in \left( \frac{T}{2}-2,\frac{T}{2}+2\right) }\int_{\left\vert
s\right\vert \leq 1/2}\left\vert W_{\ell _{o},\lambda }(x,s)\right\vert
^{2}dsd\ell _{o}dt \\
& \leq 2\int_{\ell _{o}\in \left( \frac{T}{2}-2,\frac{T}{2}+2\right)
}\int_{\left\vert s\right\vert \leq 1/2}\left\vert W_{\ell _{o},\lambda
}(x,s)\right\vert ^{2}dsd\ell _{o}\text{ .}%
\end{array}
\tag{F.10}  \label{F.10}
\end{equation}%
Consequently, from (\ref{F.8}), (\ref{F.10}) and integrating over $\Omega $,
we get the existence of $C>0$ such that%
\begin{equation}
\begin{array}{ll}
& \quad \int_{\Omega }\int_{t\in \left( \frac{T}{2}-1,\frac{T}{2}+1\right)
}\left\vert \Phi (t)u(x,t)\right\vert ^{2}dtdx \\
& \leq C\frac{1}{\lambda ^{2\beta k}}\int_{\Omega }\int_{\mathbb{R}%
}\left\vert \partial _{t}^{k}\left( \Phi \left( t\right) u(x,t)\right)
)\right\vert ^{2}dtdx \\
& \quad +4\int_{\ell _{o}\in \left( \frac{T}{2}-2,\frac{T}{2}+2\right)
}\left( \int_{\Omega }\int_{\left\vert s\right\vert \leq 1/2}\left\vert
W_{\ell _{o},\lambda }(x,s)\right\vert ^{2}dsdx\right) d\ell _{o}\text{ .}%
\end{array}
\tag{F.11}  \label{F.11}
\end{equation}

\bigskip

\noindent Now recall the following quantification result for unique
continuation of elliptic equation with Dirichlet boundary condition (Theorem
E applied to $T_{1}=-1$, $T_{2}=1$, $\delta=1/2$, $\varphi\in C_{0}^{\infty
}\left( \omega\times\left( -1,1\right) \right) $):

\bigskip

\textbf{Proposition 2 .-}\qquad \textit{Let }$\Omega $\textit{\ be a bounded
open set in }$\mathbb{R}^{n}$\textit{, }$n\geq 1$\textit{, either convex or }%
$C^{2}$\textit{\ and connected. Let }$\omega $\textit{\ be a non-empty open
subset in }$\Omega $\textit{. Let }$f=f\left( x,s\right) \in L^{2}\left(
\Omega \times \left( -1,1\right) \right) $\textit{. Then there exists }$%
\widetilde{c}>0$\textit{\ such that for all }$w=w\left( x,s\right) \in
H^{2}\left( \Omega \times \left( -1,1\right) \right) $\textit{\ solution of }%
\[
\left\{
\begin{array}{ll}
\partial _{s}^{2}w+\Delta w=f & \quad \text{\textit{in}}~\Omega \times
\left( -1,1\right) \text{ ,} \\
w=0 & \quad \text{\textit{on}}~\partial \Omega \times \left( -1,1\right)
\text{ ,}%
\end{array}%
\right.
\]%
\textit{for all }$\varepsilon >0$\textit{, we have :}%
\[
\begin{array}{ll}
& \quad \int_{\left\vert s\right\vert \leq 1/2}\int_{\Omega }\left\vert
w\left( x,s\right) \right\vert ^{2}dxds \\
& \leq \widetilde{c}e^{\widetilde{c}/\varepsilon }\left( \int_{\left\vert
s\right\vert \leq 1}\int_{\omega }\left\vert w\left( x,s\right) \right\vert
^{2}dxds+\int_{\left\vert s\right\vert \leq 1}\int_{\Omega }\left\vert
f\left( x,s\right) \right\vert ^{2}dxds\right) \\
& \quad +e^{-4c_{0}/\varepsilon }\int_{\left\vert s\right\vert \leq
1}\int_{\Omega }\left\vert w\left( x,s\right) \right\vert ^{2}dxds\text{ .}%
\end{array}%
\]

\bigskip

\noindent Applying to $W_{\ell _{0},\lambda }$, from (\ref{F.5}) we deduce
that for all $\varepsilon >0$,
\begin{equation}
\begin{array}{ll}
& \quad \int_{\left\vert s\right\vert \leq 1/2}\int_{\Omega }\left\vert
W_{\ell _{0},\lambda }(x,s)\right\vert ^{2}dxds \\
& \leq e^{-4c_{0}/\varepsilon }\int_{\left\vert s\right\vert \leq
1}\int_{\Omega }\left\vert W_{\ell _{0},\lambda }(x,s)\right\vert ^{2}dxds
\\
& \quad +\widetilde{c}e^{\widetilde{c}/\varepsilon }\int_{\left\vert
s\right\vert \leq 1}\int_{\omega }\left\vert W_{\ell _{0},\lambda
}(x,s)\right\vert ^{2}dxds \\
& \quad +\widetilde{c}e^{\widetilde{c}/\varepsilon }\int_{\left\vert
s\right\vert \leq 1}\int_{\Omega }\left\vert \int_{\mathbb{R}}-F_{\lambda
}(\ell _{0}+is-\ell )\right. \\
& \quad \quad \quad \quad \quad \quad \left. \left. \left[ \Phi ^{\prime
\prime }(\ell )u(x,\ell )+2\Phi ^{\prime }(\ell )\partial _{t}u(x,\ell )%
\right] d\ell \right\vert ^{2}dxds\right\vert ^{2}dxds\text{ .}%
\end{array}
\tag{F.12}  \label{F.12}
\end{equation}%
Consequently, from (\ref{F.11}) and (\ref{F.12}), there exists a constant $%
C>0$, such that for all $\varepsilon >0$,
\begin{equation}
\begin{array}{ll}
& \quad \int_{\Omega }\int_{t\in \left( \frac{T}{2}-1,\frac{T}{2}+1\right)
}\left\vert \Phi (t)u(x,t)\right\vert ^{2}dtdx \\
& \leq C\frac{1}{\lambda ^{2\beta k}}\int_{\Omega }\int_{\mathbb{R}%
}\left\vert \partial _{t}^{k}\left( \Phi \left( t\right) u(x,t)\right)
)\right\vert ^{2}dtdx \\
& \quad +4e^{-4c_{0}/\varepsilon }\int_{\ell _{o}\in \left( \frac{T}{2}-2,%
\frac{T}{2}+2\right) }\left( \int_{\left\vert s\right\vert \leq
1}\int_{\Omega }\left\vert W_{\ell _{0},\lambda }(x,s)\right\vert
^{2}dxds\right) d\ell _{o} \\
& \quad +4Ce^{C/\varepsilon }\int_{\ell _{o}\in \left( \frac{T}{2}-2,\frac{T%
}{2}+2\right) }\left( \int_{\left\vert s\right\vert \leq 1}\int_{\omega
}\left\vert W_{\ell _{0},\lambda }(x,s)\right\vert ^{2}dxds\right) d\ell _{o}
\\
& \quad +4Ce^{\widetilde{c}/\varepsilon }\int_{\ell _{o}\in \left( \frac{T}{2%
}-2,\frac{T}{2}+2\right) }\left( \int_{\left\vert s\right\vert \leq
1}\int_{\Omega }\left\vert \int_{\mathbb{R}}-F_{\lambda }(\ell _{0}+is-\ell
)\right. \right. \\
& \quad \quad \quad \quad \quad \quad \left. \left. \left[ \Phi ^{\prime
\prime }(\ell )u(x,\ell )+2\Phi ^{\prime }(\ell )\partial _{t}u(x,\ell )%
\right] d\ell \right\vert ^{2}dxds\right) d\ell _{o}\text{ .}%
\end{array}
\tag{F.13}  \label{F.13}
\end{equation}%
Let define $\Phi \in C_{0}^{\infty }(\mathbb{R})$ more precisely now: we
choose $\Phi \in C_{0}^{\infty }(\left( 0,T\right) )$, $0\leq \Phi \leq 1$, $%
\Phi \equiv 1$ on $\left( \frac{T}{4},\frac{3T}{4}\right) $. Furthermore,
let $K=\left[ 0,\frac{T}{4}\right] \cup \left[ \frac{3T}{4},T\right] $ such
that supp$\left( \Phi ^{\prime }\right) =K$ and supp$\left( \Phi ^{\prime
\prime }\right) \subset K$.

\bigskip

\noindent Let $K_{0}=\left[ \frac{3T}{8},\frac{5T}{8}\right] $. In
particular, dist$\left( K,K_{o}\right) =\frac{T}{8}$. Let define $T>0$ more
precisely now: we choose $T>16\max\left( 1,1/c_{2}\right) $ in order that $%
\left( \frac{T}{2}-2,\frac{T}{2}+2\right) \subset K_{0}$ and dist$\left(
K,K_{o}\right) \geq\frac{2}{c_{2}}$.

\bigskip

\noindent Now, we will choose $\ell_{0}\in\left( \frac{T}{2}-2,\frac{T}{2}%
+2\right) \subset K_{0}$ and $s\in\left[ -1,1\right] $. Consequently, for
any $\ell\in K$, $\left\vert \ell_{0}-\ell\right\vert \geq\frac{2}{c_{2}}\geq%
\frac{1}{c_{2}}\left\vert s\right\vert $ and it will imply from the second
line of (\ref{F.3}) that
\begin{equation}
\forall\ell\in K\quad\left\vert F_{\lambda}(\ell_{o}+is-\ell)\right\vert
\leq A\lambda^{\gamma}e^{-c_{1}\lambda\left( \frac{T}{8}\right) ^{1/\gamma}}%
\text{ .}  \tag{F.14}  \label{F.14}
\end{equation}

\bigskip

\noindent Till the end of the proof, $C$ and respectively $C_{T}$ will
denote a generic positive constant independent of $\varepsilon$ and $\lambda$
but dependent on $\Omega$ and respectively $\left( \Omega,T\right) $, whose
value may change from line to line.

\bigskip

\noindent The first term in the right hand side of (\ref{F.13}) becomes,
using (\ref{F.1}),
\begin{equation}
\frac{1}{\lambda ^{2\beta k}}\int_{\Omega }\int_{\mathbb{R}}\left\vert
\partial _{t}^{k}\left( \Phi \left( t\right) u(x,t)\right) )\right\vert
^{2}dtdx\leq C_{T}\frac{1}{\lambda ^{2\beta k}}\left\Vert \left(
u_{0},u_{1}\right) \right\Vert _{D\left( A^{k-1}\right) }^{2}\text{ .}
\tag{F.15}  \label{F.15}
\end{equation}%
The second term in the right hand side of (\ref{F.13}) becomes, using the
first line of (\ref{F.3}),%
\begin{equation}
\begin{array}{ll}
& \quad e^{-4/\varepsilon }\int_{\ell _{o}\in \left( \frac{T}{2}-2,\frac{T}{2%
}+2\right) }\left( \int_{\left\vert s\right\vert \leq 1}\int_{\Omega
}\left\vert W_{\ell _{0},\lambda }(x,s)\right\vert ^{2}dxds\right) d\ell _{o}
\\
& \leq \left( C_{o}\lambda ^{\gamma }e^{\lambda c_{0}}\right)
^{2}e^{-4c_{0}/\varepsilon }\int_{\ell _{o}\in \left( \frac{T}{2}-2,\frac{T}{%
2}+2\right) }\left[ \int_{\left\vert s\right\vert \leq 1}\int_{\Omega
}\left\vert \int_{0}^{T}\left\vert u(x,\ell )\right\vert d\ell \right\vert
^{2}dxds\right] \\
& \leq C_{T}\lambda ^{2\gamma }e^{2\lambda c_{0}}e^{-4c_{0}/\varepsilon
}\left\Vert \left( u_{0},u_{1}\right) \right\Vert _{H_{0}^{1}\left( \Omega
\right) \times L^{2}\left( \Omega \right) }^{2}\text{ .}%
\end{array}
\tag{F.16}  \label{F.16}
\end{equation}%
The third term in the right hand side of (\ref{F.13}) becomes, using the
first line of (\ref{F.3}),%
\begin{equation}
\begin{array}{ll}
& \quad e^{C/\varepsilon }\int_{\ell _{o}\in \left( \frac{T}{2}-2,\frac{T}{2}%
+2\right) }\left( \int_{\left\vert s\right\vert \leq 1}\int_{\omega
}\left\vert W_{\ell _{0},\lambda }(x,s)\right\vert ^{2}dxds\right) d\ell _{o}
\\
& \leq \left( C_{o}\lambda ^{\gamma }e^{\lambda c_{0}}\right)
^{2}e^{C/\varepsilon }\int_{\ell _{o}\in \left( \frac{T}{2}-2,\frac{T}{2}%
+2\right) }\left[ \int_{\left\vert s\right\vert \leq 1}\int_{\omega
}\left\vert \int_{0}^{T}\left\vert u(x,\ell )\right\vert d\ell \right\vert
^{2}dxds\right] d\ell _{o} \\
& \leq C\lambda ^{2\gamma }e^{C\lambda }e^{C/\varepsilon }\int_{\omega
}\int_{0}^{T}\left\vert u(x,t)\right\vert ^{2}dtdx\text{ .}%
\end{array}
\tag{F.17}  \label{F.17}
\end{equation}%
The fourth term in the right hand side of (\ref{F.13}) becomes, using (\ref%
{F.14}) and the choice of $\Phi $,%
\begin{equation}
\begin{array}{ll}
& \quad e^{\widetilde{c}/\varepsilon }\int_{\ell _{o}\in \left( \frac{T}{2}%
-2,\frac{T}{2}+2\right) }\left( \int_{\left\vert s\right\vert \leq
1}\int_{\Omega }\left\vert \int_{\mathbb{R}}-F_{\lambda }(\ell _{0}+is-\ell
)\right. \right. \\
& \quad \quad \quad \quad \quad \quad \left. \left. \left[ \Phi ^{\prime
\prime }(\ell )u(x,\ell )+2\Phi ^{\prime }(\ell )\partial _{t}u(x,\ell )%
\right] d\ell \right\vert ^{2}dxds\right) d\ell _{o} \\
& \leq C\left( A\lambda ^{\gamma }e^{-c_{1}\lambda \left( \frac{T}{8}\right)
^{1/\gamma }}\right) ^{2}e^{\widetilde{c}/\varepsilon }\int_{\Omega
}\left\vert \int_{K}\left( \left\vert u(x,\ell )\right\vert +\left\vert
\partial _{t}u(x,\ell )\right\vert \right) d\ell \right\vert ^{2}dx \\
& \leq C\lambda ^{2\gamma }e^{-2c_{1}\lambda \left( \frac{T}{8}\right)
^{1/\gamma }}e^{\widetilde{c}/\varepsilon }\left\Vert \left(
u_{0},u_{1}\right) \right\Vert _{H_{0}^{1}\left( \Omega \right) \times
L^{2}\left( \Omega \right) }^{2}\text{ .}%
\end{array}
\tag{F.18}  \label{F.18}
\end{equation}%
We finally obtain from (\ref{F.15}), (\ref{F.16}), (\ref{F.17}), (\ref{F.18}%
) and (\ref{F.13}), that
\begin{equation}
\begin{array}{ll}
& \quad \int_{\Omega }\int_{t\in \left( \frac{T}{2}-1,\frac{T}{2}+1\right)
}\left\vert \Phi (t)u(x,t)\right\vert ^{2}dtdx \\
& \leq C_{T}\frac{1}{\lambda ^{2\beta k}}\left\Vert \left(
u_{0},u_{1}\right) \right\Vert _{D\left( A^{k-1}\right) }^{2} \\
& \quad +C_{T}\lambda ^{2\gamma }e^{2\lambda c_{0}}e^{-4c_{0}/\varepsilon
}\left\Vert \left( u_{0},u_{1}\right) \right\Vert _{H_{0}^{1}\left( \Omega
\right) \times L^{2}\left( \Omega \right) }^{2} \\
& \quad +C\lambda ^{2\gamma }e^{C\lambda }e^{C/\varepsilon }\int_{\omega
}\int_{0}^{T}\left\vert u(x,t)\right\vert ^{2}dtdx \\
& \quad +C\lambda ^{2\gamma }e^{-2c_{1}\lambda \left( \frac{T}{8}\right)
^{1/\gamma }}e^{\widetilde{c}/\varepsilon }\left\Vert \left(
u_{0},u_{1}\right) \right\Vert _{H_{0}^{1}\left( \Omega \right) \times
L^{2}\left( \Omega \right) }^{2}\text{ .}%
\end{array}
\tag{F.19}  \label{F.19}
\end{equation}%
We begin to choose $\lambda =\frac{1}{\varepsilon }$ in order that
\begin{equation}
\begin{array}{ll}
& \quad \int_{\Omega }\int_{t\in \left( \frac{T}{2}-1,\frac{T}{2}+1\right)
}\left\vert \Phi (t)u(x,t)\right\vert ^{2}dtdx \\
& \leq \varepsilon ^{2\beta k}C_{T}\left\Vert \left( u_{0},u_{1}\right)
\right\Vert _{D\left( A^{k-1}\right) }^{2} \\
& \quad +e^{-2c_{0}/\varepsilon }\frac{1}{\varepsilon ^{2\gamma }}%
C_{T}\left\Vert \left( u_{0},u_{1}\right) \right\Vert _{H_{0}^{1}\left(
\Omega \right) \times L^{2}\left( \Omega \right) }^{2} \\
& \quad +e^{C/\varepsilon }C\int_{\omega }\int_{0}^{T}\left\vert
u(x,t)\right\vert ^{2}dtdx \\
& \quad +C\frac{1}{\varepsilon ^{2\gamma }}\exp \left( \left( -2c_{1}\left(
\frac{T}{8}\right) ^{1/\gamma }+\widetilde{c}\right) \frac{1}{\varepsilon }%
\right) \left\Vert \left( u_{0},u_{1}\right) \right\Vert _{H_{0}^{1}\left(
\Omega \right) \times L^{2}\left( \Omega \right) }^{2}\text{ .}%
\end{array}
\tag{F.20}  \label{F.20}
\end{equation}%
We finally need to choose $T>16\max \left( 1,1/c_{2}\right) $ large enough
such that $\left( -2c_{1}\left( \frac{T}{8}\right) ^{1/\gamma }+\widetilde{c}%
\right) \leq -1$ that is $8\left( \frac{1+\widetilde{c}}{2c_{1}}\right)
^{\gamma }\leq T$, to deduce the existence of $C>0$ such that for any $%
\varepsilon \in \left( 0,1\right) $,%
\begin{equation}
\begin{array}{ll}
\int_{\Omega }\int_{t\in \left( \frac{T}{2}-1,\frac{T}{2}+1\right)
}\left\vert u(x,t)\right\vert ^{2}dtdx & \leq \int_{\Omega }\int_{t\in
\left( \frac{T}{2}-1,\frac{T}{2}+1\right) }\left\vert \Phi
(t)u(x,t)\right\vert ^{2}dtdx \\
& \leq C\varepsilon ^{2\beta k}\left\Vert \left( u_{0},u_{1}\right)
\right\Vert _{D\left( A^{k-1}\right) }^{2} \\
& \quad +Ce^{C/\varepsilon }\int_{\omega }\int_{0}^{T}\left\vert
u(x,t)\right\vert ^{2}dtdx\text{ .}%
\end{array}
\tag{F.21}  \label{F.21}
\end{equation}%
Now we conclude from (\ref{F.2}), that there exist a constant $c>0$ and a
time $T>0$ large enough such that for all $\varepsilon >0$ we have
\begin{equation}
\begin{array}{ll}
& \quad \left\Vert \left( u_{0},u_{1}\right) \right\Vert _{L^{2}\left(
\Omega \right) \times H^{-1}\left( \Omega \right) }^{2} \\
& \leq e^{c/\varepsilon }\int_{\omega }\int_{0}^{T}\left\vert
u(x,t)\right\vert ^{2}dtdx+\varepsilon ^{2\beta k}\left\Vert \left(
u_{0},u_{1}\right) \right\Vert _{D\left( A^{k-1}\right) }^{2}\text{ .}%
\end{array}
\tag{F.22}  \label{F.22}
\end{equation}%
Finally, we choose
\[
\varepsilon =\left( \frac{\left\Vert \left( u_{0},u_{1}\right) \right\Vert
_{L^{2}\left( \Omega \right) \times H^{-1}\left( \Omega \right) }}{%
\left\Vert \left( u_{0},u_{1}\right) \right\Vert _{D\left( A^{k-1}\right) }}%
\right) ^{1/\left( \beta k\right) }\text{ .}
\]

\bigskip

Theorem 1.1 is deduced by applying Theorem F to $\partial _{t}u$.

\bigskip

\end{document}